\documentclass[12pt]{amsart}      
\usepackage{amssymb}
\usepackage{eucal}
\usepackage{amsmath}
\usepackage{amscd}
\usepackage[all]{xy}           

\newcommand{\nc}{\newcommand}
\newcommand{\delete}[1]{}
\nc{\dfootnote}[1]{{}}          
\nc{\ffootnote}[1]{\dfootnote{#1}}
\nc{\mfootnote}[1]{\footnote{#1}} 

\nc{\mlabel}[1]{\label{#1}}  
\nc{\mcite}[1]{\cite{#1}}  
\nc{\mref}[1]{\ref{#1}}  
\nc{\mbibitem}[1]{\bibitem{#1}} 

\delete{
\nc{\mlabel}[1]{\label{#1}  
{\hfill \hspace{1cm}{\bf{{\ }\hfill(#1)}}}}
\nc{\mcite}[1]{\cite{#1}{{\bf{{\ }(#1)}}}}  
\nc{\mref}[1]{\ref{#1}{{\bf{{\ }(#1)}}}}  
\nc{\mbibitem}[1]{\bibitem[\bf #1]{#1}} 
}
\nc{\mkeep}[1]{\marginpar{{\bf #1}}} 
\nc{\newpart}[1]{{\bf NEW: } {\bf #1}} 

\newtheorem{theorem}{Theorem}[section]
\newtheorem{prop}[theorem]{Proposition}

\newtheorem{lemma}[theorem]{Lemma}
\newtheorem{coro}[theorem]{Corollary}
\newtheorem{prop-def}{Proposition-Definition}[section]

\newtheorem{exam}[theorem]{Example} 

\nc{\dfgen}{\Omega}
\nc{\dfrel}{\Lambda}
\nc{\dfgenb}{\vec{\omega}}
\nc{\dfrelb}{\vec{\lambda}}
\nc{\dfgene}{\omega}
\nc{\dfrele}{\lambda}
\nc{\dftimes}{\,{\scriptscriptstyle \blacksquare}\,}      
\nc{\dfotimes}{\,{\scriptscriptstyle \square}\,}    

\nc{\dfl}{\succ}
\nc{\dfr}{\prec}
\nc{\dfc}{\circ}
\nc{\dfb}{\bullet}
\nc{\dft}{\star}
\nc{\dfpair}[2]{\left[ {{#1}\atop{#2}}\right]}
\nc{\dfll}{{\dfpair{\dfl}{\dfl}}}
\nc{\dflr}{{\dfpair{\dfl}{\dfr}}}
\nc{\dflc}{{\dfpair{\dfl}{\dfc}}}
\nc{\dfrl}{{\dfpair{\dfr}{\dfl}}}
\nc{\dfrr}{{\dfpair{\dfr}{\dfr}}}
\nc{\dfrc}{{\dfpair{\dfr}{\dfc}}}
\nc{\dfcl}{{\dfpair{\dfc}{\dfl}}}
\nc{\dfcr}{{\dfpair{\dfc}{\dfr}}}
\nc{\dfcc}{{\dfpair{\dfc}{\dfc}}}
\nc{\dflt}{{\dfpair{\dfl}{\dft}}}
\nc{\dfrt}{{\dfpair{\dfr}{\dft}}}
\nc{\dfct}{{\dfpair{\dfc}{\dft}}}
\nc{\dftl}{{\dfpair{\dft}{\dfl}}}
\nc{\dftr}{{\dfpair{\dft}{\dfr}}}
\nc{\dftc}{{\dfpair{\dft}{\dfc}}}
\nc{\dftt}{{\dfpair{\dft}{\dft}}}
\nc{\dflb}{{\dfpair{\dfl}{\bullet}}}
\nc{\dfrb}{{\dfpair{\dfr}{\bullet}}}
\nc{\dftb}{{\dfpair{\dft}{\bullet}}}
\nc{\dfbt}{{\dfpair{\bullet}{\dft}}}
\nc{\dfbb}{{\dfpair{\bullet}{\bullet}}}
\nc{\dfbl}{{\dfpair{\bullet}{\dfl}}}
\nc{\dfbr}{{\dfpair{\bullet}{\dfr}}}
\nc{\dfop}{\odot}
\nc{\dfoa}{\dfop^{(1)}}
\nc{\dfob}{\dfop^{(2)}}
\nc{\dfoc}{\dfop^{(3)}}
\nc{\dfod}{\dfop^{(4)}}
\nc{\tstar}{\tilde{\star}}
\nc{\tdfop}{\tilde{\dfop}}
\nc{\tdfoa}{\tilde{\dfop}^{(1)}}
\nc{\tdfob}{\tilde{\dfop}^{(2)}}
\nc{\tdfoc}{\tilde{\dfop}^{(3)}}
\nc{\tdfod}{\tilde{\dfop}^{(4)}}
\nc{\dflls}{{{\dfl}\atop{\dfl}}}
\nc{\dflrs}{{{\dfl}\atop{\dfr}}}
\nc{\dflcs}{{{\dfl}\atop{\dfc}}}
\nc{\dfrls}{{{\dfr}\atop{\dfl}}}
\nc{\dfrrs}{{{\dfr}\atop{\dfr}}}
\nc{\dfrcs}{{{\dfr}\atop{\dfc}}}
\nc{\dfcls}{{{\dfc}\atop{\dfl}}}
\nc{\dfcrs}{{{\dfc}\atop{\dfr}}}
\nc{\dfccs}{{{\dfc}\atop{\dfc}}}
\nc{\dflts}{{{\dfl}\atop{\dft}}}
\nc{\dfrts}{{{\dfr}\atop{\dft}}}
\nc{\dfcts}{{{\dfc}\atop{\dft}}}
\nc{\dftls}{{{\dft}\atop{\dfl}}}
\nc{\dftrs}{{{\dft}\atop{\dfr}}}
\nc{\dftcs}{{{\dft}\atop{\dfc}}}
\nc{\dftts}{{{\dft}\atop{\dft}}}
\nc{\ndfll}{{\dfpair{\dfl}{>}}}
\nc{\ndflr}{{\dfpair{\dfl}{<}}}
\nc{\ndflc}{{\dfpair{\dfl}{\bullet}}}
\nc{\ndfrl}{{\dfpair{\dfr}{>}}}
\nc{\ndfrr}{{\dfpair{\dfr}{<}}}
\nc{\ndfrc}{{\dfpair{\dfr}{\bullet}}}
\nc{\ndfcl}{{\dfpair{\dfc}{>}}}
\nc{\ndfcr}{{\dfpair{\dfc}{<}}}
\nc{\ndfcc}{{\dfpair{\dfc}{\bullet}}}
\nc{\ndflt}{{\dfpair{\dfl}{\diamondsuit}}}
\nc{\ndfrt}{{\dfpair{\dfr}{\diamondsuit}}}
\nc{\ndfct}{{\dfpair{\dfc}{\diamondsuit}}}
\nc{\ndftl}{{\dfpair{\dft}{>}}}
\nc{\ndftr}{{\dfpair{\dft}{<}}}
\nc{\ndftc}{{\dfpair{\dft}{\bullet}}}
\nc{\ndftt}{{\dfpair{\dft}{\diamondsuit}}}
\nc{\dftri}[3]{{\left[{{#1}\atop{#2}}\atop{\scriptstyle{#3}}\right]}}
\nc{\dfrrr}{{\dftri{\dfr}{\dfr}{\dfr}}}
\nc{\dfrrl}{{\dftri{\dfr}{\dfr}{\dfl}}}
\nc{\dfrrb}{{\dftri{\dfr}{\dfr}{\dfb}}}
\nc{\dfrrt}{{\dftri{\dfr}{\dfr}{\dft}}}
\nc{\dfrlr}{{\dftri{\dfr}{\dfl}{\dfr}}}
\nc{\dfrll}{{\dftri{\dfr}{\dfl}{\dfl}}}
\nc{\dfrlb}{{\dftri{\dfr}{\dfl}{\dfb}}}
\nc{\dfrlt}{{\dftri{\dfr}{\dfl}{\dft}}}

\nc{\dfrbr}{{\dftri{\dfr}{\dfb}{\dfr}}}
\nc{\dfrbl}{{\dftri{\dfr}{\dfb}{\dfl}}}
\nc{\dfrbb}{{\dftri{\dfr}{\dfb}{\dfb}}}
\nc{\dfrbt}{{\dftri{\dfr}{\dfb}{\dft}}}
\nc{\dfrtr}{{\dftri{\dfr}{\dft}{\dfr}}}
\nc{\dfrtl}{{\dftri{\dfr}{\dft}{\dfl}}}
\nc{\dfrtb}{{\dftri{\dfr}{\dft}{\dfb}}}
\nc{\dfrtt}{{\dftri{\dfr}{\dft}{\dft}}}

\nc{\dflrr}{{\dftri{\dfl}{\dfr}{\dfr}}}
\nc{\dflrl}{{\dftri{\dfl}{\dfr}{\dfl}}}
\nc{\dflrb}{{\dftri{\dfl}{\dfr}{\dfb}}}
\nc{\dflrt}{{\dftri{\dfl}{\dfr}{\dft}}}
\nc{\dfllr}{{\dftri{\dfl}{\dfl}{\dfr}}}
\nc{\dflll}{{\dftri{\dfl}{\dfl}{\dfl}}}
\nc{\dfllb}{{\dftri{\dfl}{\dfl}{\dfb}}}
\nc{\dfllt}{{\dftri{\dfl}{\dfl}{\dft}}}

\nc{\dflbr}{{\dftri{\dfl}{\dfr}{\dfr}}}
\nc{\dflbl}{{\dftri{\dfl}{\dfr}{\dfl}}}
\nc{\dflbb}{{\dftri{\dfl}{\dfr}{\dfb}}}
\nc{\dflbt}{{\dftri{\dfl}{\dfr}{\dft}}}
\nc{\dfltr}{{\dftri{\dfl}{\dft}{\dfr}}}
\nc{\dfltl}{{\dftri{\dfl}{\dft}{\dfl}}}
\nc{\dfltb}{{\dftri{\dfl}{\dfl}{\dfb}}}
\nc{\dfltt}{{\dftri{\dfl}{\dft}{\dft}}}

\nc{\dfbrr}{{\dftri{\dfb}{\dfr}{\dfr}}}
\nc{\dfbrl}{{\dftri{\dfb}{\dfr}{\dfl}}}
\nc{\dfbrb}{{\dftri{\dfb}{\dfr}{\dfb}}}
\nc{\dfbrt}{{\dftri{\dfb}{\dfr}{\dft}}}
\nc{\dfblr}{{\dftri{\dfb}{\dfl}{\dfr}}}
\nc{\dfbll}{{\dftri{\dfb}{\dfl}{\dfl}}}
\nc{\dfblb}{{\dftri{\dfb}{\dfl}{\dfb}}}
\nc{\dfblt}{{\dftri{\dfb}{\dfl}{\dft}}}

\nc{\dfbbr}{{\dftri{\dfb}{\dfb}{\dfr}}}
\nc{\dfbbl}{{\dftri{\dfb}{\dfb}{\dfl}}}
\nc{\dfbbb}{{\dftri{\dfb}{\dfb}{\dfb}}}
\nc{\dfbbt}{{\dftri{\dfb}{\dfb}{\dft}}}
\nc{\dfbtr}{{\dftri{\dfb}{\dft}{\dfr}}}
\nc{\dfbtl}{{\dftri{\dfb}{\dft}{\dfl}}}
\nc{\dfbtb}{{\dftri{\dfb}{\dft}{\dfb}}}
\nc{\dfbtt}{{\dftri{\dfb}{\dft}{\dft}}}

\nc{\dftrr}{{\dftri{\dft}{\dfr}{\dfr}}}
\nc{\dftrl}{{\dftri{\dft}{\dfr}{\dfl}}}
\nc{\dftrb}{{\dftri{\dft}{\dfr}{\dfb}}}
\nc{\dftrt}{{\dftri{\dft}{\dfr}{\dft}}}
\nc{\dftlr}{{\dftri{\dft}{\dfl}{\dfr}}}
\nc{\dftll}{{\dftri{\dft}{\dfl}{\dfl}}}
\nc{\dftlb}{{\dftri{\dft}{\dfl}{\dfb}}}
\nc{\dftlt}{{\dftri{\dft}{\dfl}{\dft}}}

\nc{\dftbr}{{\dftri{\dft}{\dfb}{\dfr}}}
\nc{\dftbl}{{\dftri{\dft}{\dfb}{\dfl}}}
\nc{\dftbb}{{\dftri{\dft}{\dfb}{\dfb}}}
\nc{\dftbt}{{\dftri{\dft}{\dfb}{\dft}}}
\nc{\dfttr}{{\dftri{\dft}{\dft}{\dfr}}}
\nc{\dfttl}{{\dftri{\dft}{\dft}{\dfl}}}
\nc{\dfttb}{{\dftri{\dft}{\dft}{\dfb}}}
\nc{\dfttt}{{\dftri{\dft}{\dft}{\dft}}}

\newcommand{\se}{\searrow}

\newcommand{\sep}{\se^{op}}

\nc{\bin}[2]{ (_{\stackrel{\scs{#1}}{\scs{#2}}})}  
\nc{\binc}[2]{ \left (\!\! \begin{array}{c} \scs{#1}\\
    \scs{#2} \end{array}\!\! \right )}  
\nc{\bincc}[2]{  \left ( {\scs{#1} \atop
    \vspace{-.5cm}\scs{#2}} \right )}  
\nc{\bs}{\bar{S}}
\nc{\la}{\longrightarrow}
\nc{\rar}{\rightarrow}
\nc{\dar}{\downarrow}
\nc{\dap}[1]{\downarrow \rlap{$\scriptstyle{#1}$}}
\nc{\uap}[1]{\uparrow \rlap{$\scriptstyle{#1}$}}
\nc{\defeq}{\stackrel{\rm def}{=}}
\nc{\oeq}[1]{\stackrel{(#1)}{=}}
\nc{\dis}[1]{\displaystyle{#1}}
\nc{\dotcup}{\ \displaystyle{\bigcup^\bullet}\ }
\nc{\hcm}{\ \hat{,}\ }
\nc{\hcirc}{\hat{\circ}}
\nc{\hts}{\hat{\otimes}}
\nc{\lts}{\stackrel{\leftarrow}{\otimes}}
\nc{\rts}{\stackrel{\rightarrow}{\otimes}}
\nc{\lleft}{[}
\nc{\lright}{]}
\nc{\curlyl}{\left \{ \begin{array}{c} {} \\ {} \end{array}
    \right .  \!\!\!\!\!\!\!}
\nc{\curlyr}{ \!\!\!\!\!\!\!
    \left . \begin{array}{c} {} \\ {} \end{array}
    \right \} }
\nc{\longmid}{\left | \begin{array}{c} {} \\ {} \end{array}
    \right . \!\!\!\!\!\!\!}
\nc{\ora}[1]{\stackrel{#1}{\rar}}
\nc{\ola}[1]{\stackrel{#1}{\la}}
\nc{\olrarrow}[1]{\stackrel{#1}{\Leftrightarrow}}
\nc{\scs}[1]{\scriptstyle{#1}}
\nc{\mrm}[1]{{\rm #1}}
\nc{\margin}[1]{\marginpar{\rm #1}}   
\nc{\dirlim}{\displaystyle{\lim_{\longrightarrow}}\,}
\nc{\invlim}{\displaystyle{\lim_{\longleftarrow}}\,}
\nc{\mvp}{\vspace{0.5cm}}
\nc{\ot}{\otimes}
\nc{\otl}{\otimes_\ell}
\nc{\otr}{\otimes_r}
\nc{\tk}{^{(k)}}
\nc{\tp}{^\prime}
\nc{\ttp}{^{\prime\prime}}
\nc{\svp}{\vspace{2cm}}
\nc{\vp}{\vspace{8cm}}
\nc{\proofbegin}{\noindent{\bf Proof: }}
\nc{\proofend}{$\blacksquare$ \vspace{0.5cm}}
\nc{\modg}[1]{\!<\!\!{#1}\!\!>}
\nc{\intg}[1]{F_C(#1)}
\nc{\lmodg}{\!<\!\!}
\nc{\rmodg}{\!\!>\!}
\nc{\cpi}{\widehat{\Pi}}
\nc{\sha}{{\mbox{\cyr X}}}  
\nc{\shpr}{\diamond}    
\nc{\vep}{\varepsilon}
\nc{\labs}{\mid\!}
\nc{\rabs}{\!\mid}
\nc{\hsha}{\widehat{\sha}}
\nc{\lsha}{\stackrel{\leftarrow}{\sha}}
\nc{\rsha}{\stackrel{\rightarrow}{\sha}}

\nc{\ann}{\mrm{ann}}
\nc{\Aut}{\mrm{Aut}}
\nc{\can}{\mrm{can}}
\nc{\coh}{\mrm{coh}}
\nc{\comp}{\mrm{comp}}
\nc{\colim}{\mrm{colim}}
\nc{\Cont}{\mrm{Cont}}
\nc{\rchar}{\mrm{char}}
\nc{\cok}{\mrm{coker}}
\nc{\dtf}{{R-{\rm tf}}}
\nc{\dtor}{{R-{\rm tor}}}

\nc{\Div}{{\mrm Div}}
\nc{\End}{\mrm{End}}
\nc{\Ext}{\mrm{Ext}}
\nc{\Fil}{\mrm{Fil}}
\nc{\Frob}{\mrm{Frob}}
\nc{\Gal}{\mrm{Gal}}
\nc{\GL}{\mrm{GL}}
\nc{\Hom}{\mrm{Hom}}
\nc{\hsr}{\mrm{H}}
\nc{\hpol}{\mrm{HP}}
\nc{\id}{\mrm{id}}
\nc{\im}{\mrm{im}}
\nc{\incl}{\mrm{incl}}
\nc{\length}{\mrm{length}}
\nc{\Loday}{ABQR\ }
\nc{\mchar}{\rm char}
\nc{\mpart}{\mrm{part}}
\nc{\ql}{{\QQ_\ell}}
\nc{\qp}{{\QQ_p}}
\nc{\rank}{\mrm{rank}}
\nc{\rcot}{\mrm{cot}}
\nc{\rdef}{\mrm{def}}
\nc{\rdiv}{{\rm div}}
\nc{\rtf}{{\rm tf}}
\nc{\rtor}{{\rm tor}}
\nc{\res}{\mrm{res}}
\nc{\SL}{\mrm{SL}}
\nc{\Spec}{\mrm{Spec}}
\nc{\tor}{\mrm{tor}}
\nc{\Tr}{\mrm{Tr}}
\nc{\tr}{\mrm{tr}}
\nc{\vect}{\mathbf{Vect}}

\nc{\ab}{\mathbf{Ab}}
\nc{\Alg}{\mathbf{Alg}}
\nc{\Bax}{\mathbf{Bax}}
\nc{\bfk}{{\bf k}}
\nc{\bfone}{{\bf 1}}
\nc{\base}[1]{{a_{#1}}}
\nc{\detail}{\marginpar{\bf More detail}
    \noindent{\bf Need more detail!}
    \svp}
\nc{\Diff}{\mathbf{Diff}}
\nc{\gap}{\marginpar{\bf Incomplete}\noindent{\bf Incomplete!!}
    \svp}
\nc{\FMod}{\mathbf{FMod}}
\nc{\Int}{\mathbf{Int}}
\nc{\Mon}{\mathbf{Mon}}
\nc{\remarks}{\noindent{\bf Remarks: }}
\nc{\Rep}{\mathbf{Rep}}
\nc{\Rings}{\mathbf{Rings}}
\nc{\Sets}{\mathbf{Sets}}

\nc{\BA}{{\Bbb A}}
\nc{\CC}{{\Bbb C}}
\nc{\DD}{{\Bbb D}}
\nc{\EE}{{\Bbb E}}
\nc{\FF}{{\Bbb F}}
\nc{\GG}{{\Bbb G}}
\nc{\HH}{{\Bbb H}}
\nc{\LL}{{\Bbb L}}
\nc{\NN}{{\Bbb N}}
\nc{\QQ}{{\Bbb Q}}
\nc{\RR}{{\Bbb R}}
\nc{\TT}{{\Bbb T}}
\nc{\VV}{{\Bbb V}}
\nc{\ZZ}{{\Bbb Z}}

\nc{\cala}{{\mathcal A}}
\nc{\calc}{{\mathcal C}}
\nc{\cald}{{\mathcal D}}
\nc{\cale}{{\mathcal E}}
\nc{\calf}{{\mathcal F}}
\nc{\calg}{{\mathcal G}}
\nc{\calh}{{\mathcal H}}
\nc{\cali}{{\mathcal I}}
\nc{\call}{{\mathcal L}}
\nc{\calm}{{\mathcal M}}
\nc{\caln}{{\mathcal N}}
\nc{\calo}{{\mathcal O}}
\nc{\calp}{{\mathcal P}}
\nc{\calq}{{\mathcal Q}}
\nc{\calr}{{\mathcal R}}
\nc{\calt}{{\mathcal T}}
\nc{\calw}{{\mathcal W}}
\nc{\calx}{{\mathcal X}}
\nc{\CA}{\mathcal{A}}

\nc{\fraka}{{\frak a}}
\nc{\frakB}{{\frak B}}
\nc{\frakm}{{\frak m}}
\nc{\frakp}{{\frak p}}

\font\cyr=wncyr10

\begin{document}

\title{Coherent Unit Actions on Regular Operads and Hopf Algebras
}
\author{Kurusch Ebrahimi-Fard}
\address{Institut Henri Poincare,
         11, rue Pierre et Marie Curie
         F-75231 Paris Cedex 05
         FRANCE
and      Universit\"at Bonn -
         Physikalisches Institut,
         Nussallee 12,
         D-53115 Bonn,
         Germany}
\email{kurusch@ihes.fr}
\author{Li Guo}
\address{
Department of Mathematics and Computer Science,
Rutgers University,
Newark, NJ 07102, USA}
\email{liguo@newark.rutgers.edu}

\date{\today}

\begin{abstract}
J.-L. Loday introduced the concept of coherent unit actions on a
regular operad and showed that such actions give
Hopf algebra structures on the free algebras. Hopf algebras
obtained this way include the Hopf algebras of shuffles,
quasi-shuffles and planar rooted trees. We characterize coherent
unit actions on binary quadratic regular operads in terms of
linear equations of the generators of the operads. We then use
these equations to classify operads with coherent unit actions. We
further show that coherent unit actions are preserved under taking
products and thus yield Hopf algebras on the free object of the
product operads when the factor operads have coherent unit
actions. On the other hand, coherent unit actions are never
preserved under taking the dual in the operadic sense except for
the operad of associative algebras.
\end{abstract}

\delete{
\begin{keyword}
coherent unit actions\sep regular operads\sep Hopf algebras\sep dendriform algebras.

\PACS 18D50 \sep 17A30 \sep 16W30

\end{keyword}
}

\maketitle


\setcounter{section}{0}

\section{Introduction}
\label{intro}

While the original motivation for the study of the
dendriform dialgebra~\cite{Lo1,Lo2} was to study the
periodicity of algebraic $K$-groups, it soon became clear that
dendriform dialgebras are an interesting subject on its own. This
can be seen on one hand by its quite extensive study by several
authors in areas related to operads~\cite{Lo6},
homology~\cite{Fra1,Fra2}, combinatorics~\cite{A-S1,A-S2,Fo,L-R1},
arithmetic~\cite{Lo5}, quantum field theory~\cite{Fo} and
especially Hopf algebras~\cite{A-S1,Ch,Hol1,L-R3,Ron}. On the
other hand it has several generalizations and extensions that
share properties of the original dendriform dialgebra. These
new structures include the dendriform trialgebra~\cite{L-R2},
the dipterous algebra~\cite{L-R3}, the
dendriform quadri-algebra~\cite{A-L}, the 2-associative
algebra~\cite{L-R3,L-R4,Pi}, the magma algebra~\cite{G-H} and
the ennea-algebra. In fact, they are special cases of a class of
binary quadratic regular operads that will be made precise later
in this paper.

It is remarkable that many of these algebra structures have a Hopf
algebra structure on the free algebras. For example, the free
commutative dendriform dialgebras and trialgebras are the shuffle
and quasi-shuffle Hopf algebras, and the free dendriform dialgebra
and trialgebras are the Hopf algebra of binary planar rooted
trees~\mcite{Lo4,L-R1}  and planar rooted trees~\mcite{L-R2}.
These findings were put in a general framework recently by
Loday~\mcite{Lo6} who showed that the existence of a coherent unit
action on a binary quadratic regular operad with a splitting of
associativity endows the free objects with a Hopf algebra
structure. Since then, this method has been applied to obtain Hopf
algebra structures on several other operads~\cite{Le1,Le2,Le3,Lo6}.

Thus it is desirable to obtain a good understanding of such operads with
a coherent unit action and therefore with a Hopf algebra structure
on the free objects. This is our goal to achieve in this paper,
by working
with the generators and relations of these regular operads~\cite{E-G1,Lo7}.
As a result, we explicitly describe a large
class of operads that give rise to Hopf algebras.

After briefly recalling related concepts and results, we first
give in Section~\ref{sec:def} a simple criterion for a unit action
to be coherent. This criterion reduces the checking of the coherence
condition to the verification of a system of linear equations,
called coherence equations. Then in Section~\ref{sec:coh} we use
the coherence equations to obtain a classification of binary, quadratic,
regular operads that allows a coherent unit action.
Special cases are studied and are
related to examples in the current literature.

The compatibility of coherent unit actions on operads with taking
operad products and duals is studied in
Section~\ref{sec:prod}. We show that the coherence condition is
preserved by taking products. Thus the Hopf algebra structure on
the product operad follows automatically from those on the factor
operads, as long as the factor operads have coherent unit actions.
In contrast to products, we
show that the coherence condition is never preserved by taking the
dual in the operadic sense, except for the trivial case when the
operad is the one for associative algebras.

We give a similar study of the related notion of compatible unit actions
~\cite{Lo6}. It is related to, but weaker than the notion of coherent
unit actions.

\section{Compatible and coherent unit actions}
\mlabel{sec:def}

\subsection{\Loday operads}
We recall the standard definition of algebraic operads in general
before rephrasing it in our special case.
Since we will not need the general definition in the rest of the
paper, we refer the interested reader to find further details in
the standard references, such as~\cite{G-K,Lo3,L-S-V,M-S-S}.

Let $\bfk$ be a field of characteristic zero and let $\vect$ be
the category of $\bfk$-vector spaces. An algebraic operad over
$\bfk$ is an analytic functor $\calp:\vect\to \vect$ such that
$\calp(0)=0$, and is equipped with a natural transformation of
functors $\gamma: \calp\circ \calp \to \calp$ which is associative
and has a unit $1:\id \to \calp$.

By considering free $\calp$-algebras, an operad gives a sequence $\{\calp(n)\}$ of
finitely generated $\bfk[S_n]$-modules that satisfy certain composition axioms.
An operad is called {\bf binary} if $\calp(1)=\bfk$ and $\calp(2)$ generates
$\calp(n), n\geq 3$ by composition; is called {\bf quadratic} if
all relations among the binary operations $\calp(2)$ are derived from
relations in $\calp(3)$; is called {\bf regular} if,
moreover, the binary operations have no symmetries and the
variables $x,\, y$ and $z$ appear in the same order
(such as $(x\cdot y) \cdot z=x\cdot
(y\cdot z)$, not $(x\cdot y)\cdot z=x\cdot (z\cdot y)$).

By regularity, the space $\calp(n)$ is of the form $\calp_n\otimes
\bfk[S_n]$ where $\calp_n$ is a vector space. So the operad
$\{\calp(n)\}$ is determined by $\{\calp_n\}$. Then a binary,
quadratic, regular operad is determined by a pair
$(\dfgen,\dfrel)$ where $\dfgen=\calp_2$, called the {\bf space of
generators}, and $\dfrel$ is a subspace of $\dfgen^{\otimes 2}
\oplus \dfgen^{\otimes 2}$, called the {\bf space of relations}.
So we write $\calp=(\dfgen,\dfrel)$ for the operad.
Since $(\dfgen,\dfrel)$ is determined by
$(\dfgenb,\dfrelb)$ where $\dfgenb$ is a basis of $\dfgen$ and
$\dfrelb$ is a basis of $\dfrel$, we also use $(\dfgenb,\dfrelb)$
to denote a binary, quadratic, regular operad, as is
usually the case in the literature.

For such a $\calp=(\dfgen,\dfrel)$, a $\bfk$-vector space $A$ is called
a $\calp$-algebra if it has binary operations $\dfgen$ and if,
for
$$\big (\sum_{i=1}^k\dfoa_i\otimes \dfob_i, \sum_{j=1}^m\dfoc_j\otimes \dfod_j \big)
\in \dfrel \subseteq
\dfgen^{\otimes 2} \oplus \dfgen^{\otimes 2}$$ with
$\dfoa_i,\dfob_i,\dfoc_j,\dfod_j\in \dfgen$, $1\leq i\leq k$, $1\leq j\leq m$,
we have
\begin{equation}
 \sum_{i=1}^k(x\dfoa_i y) \dfob_i z = \sum_{j=1}^m x \dfoc_j (y \dfod_j z), \forall\ x,y,z\in A.
\mlabel{eq:rel}
\end{equation}

We say that a binary, quadratic, regular operad $(\dfgen,\dfrel)$
{\bf has a splitting associativity} if there is a choice of
$\star$
in $\dfgen$ such that $(\star\otimes \star, \star\otimes \star)$
is in $\dfrel$~\cite{Lo6}. As an abbreviation, we call such an operad
an {\bf associative BQR operad}, or simply an {\bf \Loday operad}.
Equivalently~\cite[Lemma 2.1]{E-G1}, a binary, quadratic, regular operad
is \Loday if and
only if there is a basis $\dfgenb=\{\dfgene_i\}$ of $\dfgen$ such
that $\star=\sum_i \dfgene_i$ and there is a basis
$\dfrelb=\{\dfrele_j\}_j$ of $\dfrel$ such that the associativity
of $\star$ is given by the sum of $\dfrele_j$, giving a splitting of the
associativity of $\star$\,:
$$ (\star\otimes \star, \star\otimes \star)= \sum_j \dfrele_j.$$
Note that a binary quadratic regular operad $(\dfgen,\dfrel)$
might have different choices of associative operations.
For example, if $\star$ is associative,
then so is $c\star$ for any nonzero $c\in\bfk$. As we will see
later, some operads even have linear independent associative
operations.

To be precise, we let $(\dfgen,\dfrel,\star)$, or $(\dfgenb,\dfrelb,\star)$,
denote an \Loday operad with
$\star$ as the chosen associative operation.
Let $(\dfgen,\dfrel,\star)$ and $(\dfgen',\dfrel',\star')$ be \Loday operads with
associative operations $\star$ and $\star'$ respectively.
A {\bf morphism} $f:(\dfgen,\dfrel,\star)\to (\dfgen',\dfrel',\star)$ is
a linear map from $\dfgen$ to $\dfgen'$ sending $\star$ to $\star'$ and
inducing a linear map from $\dfrel$ to $\dfrel'$. An invertible morphism
is called an {\bf isomorphism}.

The following examples of \Loday operads will be used later in the paper.
\begin{exam}{\rm
\begin{enumerate}
\item An associative $\bfk$-algebra
is a $\bfk$-vector space $A$ with an associative \mbox{product
$\cdot$\,.} The corresponding operad is
$\calp_A=(\dfgenb_A,\dfrelb_A, \cdot)$
with $\dfgenb_A=\{\cdot\}$ and
$\dfrelb_A=\{ (\cdot\otimes\cdot,\cdot\otimes\cdot)\}.$
\item The {\bf dendriform dialgebra} of Loday~\cite{Lo4} corresponds
to the operad $\calp_D=(\dfgenb_D,\dfrelb_D,\star_D)$ with
$\dfgenb_D=\{\prec,\succ\}$, $\star_D=\prec+ \succ$ and
\begin{equation}
\dfrelb_D=\{(\prec\otimes\prec,\prec\otimes(\prec+\succ)),
(\succ\otimes\prec,\succ\otimes\prec),
((\prec+\succ)\otimes\succ,\succ\otimes\succ)\}.
\mlabel{eq:dia}
\end{equation}
\item The {\bf dendriform trialgebra} of Loday and
Ronco~\cite{L-R2} corresponds to the operad
$\calp_T=(\dfgenb_T,\dfrelb_T, \star_T)$ with
$\dfgenb_T=\{\prec,\succ,\circ\}$, $\star_T=\prec+\succ+\,\circ$ and
\begin{eqnarray}
\dfrelb_T=& \{(\prec\otimes\prec,\prec\otimes\:\star),
(\succ\otimes\prec,\succ\otimes\prec), (\star\:\otimes
\succ,\succ\otimes\succ),
(\succ\otimes\:\circ,\succ\otimes\:\circ), \notag\\
&
(\prec\otimes\:\circ,\circ\:\otimes\succ),(\circ\:\otimes\prec,\circ\:\otimes\prec),
(\circ\otimes\circ,\circ\otimes\circ)\}. \mlabel{eq:tri}
\end{eqnarray}
\item Leroux's {\bf NS-algebra} \cite{Le2}
corresponds to the operad $\calp_N=(\dfgenb_N, \dfrelb_N,\star_N)$ with
$\dfgenb_N=\{\prec,\succ ,\bullet\}$, $\star_N=\prec+\succ+\,\bullet$ and
\begin{eqnarray}
\dfrelb_N=& \{(\prec\otimes\prec,\prec\otimes\:\star),
(\succ\otimes\prec,\succ\otimes\prec),
(\star\:\otimes \succ,\succ\otimes\succ), \nonumber\\
& (\star\:\otimes\bullet+\bullet\:\otimes\prec,
\succ\otimes\bullet+\bullet\otimes\:\star)\}. \mlabel{eq:NS}
\end{eqnarray}

\end{enumerate}
}
%
%
\end{exam}

\subsection{Compatible and coherent unit actions}
\mlabel{sec:unit} We now review the concepts of coherent and
compatible actions on regular operads, and the theorem of Loday
showing that the existence of coherent
unit actions yields Hopf algebras.
\medskip

Let $\calp=(\dfgen,\dfrel,\star)$ be an \Loday operad. A {\bf unit action} on $\calp$ is
a choice of two linear maps
$$\alpha,\ \beta: \dfgen \xrightarrow{} \bfk$$
such that
$\alpha(\star)=\beta(\star)=1,$ the unit of $\bfk$.

Let $A$ be a $\calp$-algebra. A unit action $(\alpha,\beta)$ allows us to
extend a binary operation $\dfop\in \dfgen$ on $A$ to a restricted binary operation on
$A_+:=\bfk . 1\oplus A$ by defining
\begin{equation}
\tdfop: A_+ \otimes A_+ \to A_+, \ \
a\tdfop b:= \left \{\begin{array}{ll} a\dfop b, & a,b\in A,\\
\alpha(\dfop)\, a, & a\in A, b=1, \\
\beta(\dfop)\, b, & a=1, b\in A, \\
1, & a=b=1, \dfop =\star, \\
{\rm undefined}, & a=b=1, \dfop\neq \star.
\end{array} \right .
\mlabel{eq:ext}
\end{equation}
Thus the extended binary operation $\tdfop$ is defined on the
subspace $(\bfk. 1\ot A)\oplus (A\ot \bfk. 1) \oplus (A\ot A)$ of
$A_+\ot A_+$, and on the full space $A_+\otimes A_+$ when $\dfop
=\star$. The unit action is called {\bf compatible}
if the relations of $\calp$ are still valid
on $A_+$ for each $\calp$-algebra $A$ whenever the terms are defined.
More precisely,
if $(\sum_{i} \dfoa_i\ot \dfob_i, \sum_j \dfoc_j\ot \dfod_j)$ is a relation
of $\calp$, then
$$ \sum_i (x\tilde{\dfop}^{(1)}_i y) \tilde{\dfop}^{(2)}_i z
    =\sum_j x\tilde{\dfop}^{(3)}_j (y \tilde{\dfop}^{(4)}_j z),\
    \forall\, x,y,z\in A_+,$$
    whenever the two sides make sense.

Next let $A$, $B$ be two $\calp$-algebras where $\calp$ has a compatible
unit action $(\alpha,\beta)$. Consider the subspace
$A\boxtimes B:= (A \otimes \bfk.1) \oplus (\bfk.1 \otimes B) \oplus (A \otimes B)$.
For $\dfop\in \dfgen$,
define a binary operation
$$\boxdot: (A\boxtimes B) \otimes (A\boxtimes B)\to A\boxtimes B$$
by
\begin{equation}
(a \otimes b) \boxdot (a' \otimes b') :=
\left \{ \begin{array}{ll} (a \tstar a') \otimes (b \tdfop b'),  &
    \textrm{if} \ \ b \otimes b' \not= 1 \otimes 1, \\
(a \tdfop a') \otimes 1, &  \textrm{otherwise}.
\end{array} \right .
\mlabel{eq:tensorext}
\end{equation}
We say that the unit action $(\alpha,\beta)$ is {\bf{coherent}}
if, for any $\calp$-algebras $A$ and $B$, the subspace
 $A \boxtimes B$, equipped with these operations is still a $\calp$-algebra.
We note that $A_+\ot B_+=\bfk.1 \oplus (A\boxtimes B)=(A\boxtimes B)_+$. Thus
the associative operation $\star$ gives an associative operation
$\star$ on $A\boxtimes B$ which,
as in (\ref{eq:ext}), extends to an associative operation $\tstar$
on $A_+\ot B_+$.

By definition, a coherent unit action on $\calp$ is also compatible.
The converse is not true. See Example~\mref{ex:twoass}.

The significance of the coherence property can be seen in the following theorem
of Loday~\cite{Lo6}. We refer the reader to the original paper for further details.
\begin{theorem}
Let $\calp$ be an \Loday operad. Let $\calp(V)_+$ be the augmented free $\calp$-algebra
on a $\bfk$-vector space $V$.
Any coherent unit action on $\calp$ equips $\calp(V)_+$
with a connected Hopf algebra structure.
\label{thm:Loday}
\end{theorem}
The Hopf algebra structure is in fact a $\calp$-Hopf algebra structure in the
sense that the coproduct is a morphism of augmented $\calp$-algebras.

For example, it was shown in~\cite{Lo6} that,
for the dendriform trialgebra
$\calp_T=(\dfgenb_T,\dfrelb_T,\star_T)$ with $\dfgenb_T=\{\prec,
\succ,\circ\}$, the unit action
$$ \alpha(\prec)=\beta(\succ)=1,
\alpha(\succ)=\alpha(\circ)=\beta(\prec)=\beta(\circ)=0$$ is
coherent with the relations of $\calp_T$, and thus the free
$\calp_T$-algebra on a $\bfk$-vector space has a Hopf algebra
structure. This is the Hopf algebra of planar rooted
trees~\mcite{L-R1}. The same method also recovers the Hopf algebra
structure on the free dendriform dialgebra as planar binary rooted
trees~\cite{Lo4}, and applies to obtain Hopf algebra structures on
several other algebras~\cite{Le1,Le2,Le3,Lo6}
(see
Corollary~\ref{co:Hopf}). Furthermore, a variation of
Theorem~\ref{thm:Loday} applies to Zinbiel algebras to recover the
shuffle Hopf algebra~\cite{Lo6} and applies to commutative
trialgebras to recover the quasi-shuffle Hopf algebra of
Hoffman~\mcite{Ho}. See~\mcite{Lo8}.

The understandings gained in this paper on coherent unit actions on
\Loday operads, especially the classification (Theorem~\mref{thm:crel}),
give us the precise information on the kind of operads that we should expect
a Hopf algebra structure on the free objects.
For some of these operads, the construction of their free objects might be
obtained through the relation between the operads and certain type of
linear operators. For example, there are natural functors from the category
of Rota-Baxter algebras (of weight zero and one) to
that of dendriform dialgebras and
trialgebras~\mcite{Ag1,EF1}. By an analog of the Poincar\'e-Birkhoff-Witt
Theorem~\mcite{E-G2}, the adjoint functors embed a dendriform
dialgebra and trialgebra, especially a free one, into a Rota-Baxter
algebra. Thus one can use Rota-Baxter algebras to
construct free dendriform dialgebras and trialgebras and recover the multiplication
in the Hope algebra structures. A Similar approach should apply to some other instances,
such as the natural functor from Nijenhuis algebras to Lerous'
NS-algebras~\mcite{Le2}
and the natural functors from algebras with two commuting Rota-Baxter
operators to dendriform quadri-algebras and ennea-algebras~\mcite{A-L,Le1}.

Coherent unit actions are also defined by Loday~\cite{Lo6} for binary,
quadratic operads without the regularity condition
and are recently extended to
general algebraic operads by Holtkamp~\cite{Hol2}. In these cases,
the free objects of the operad $\calp$ have the structure of a
$\calp$-Hopf algebra, a more general concept than Hopf algebra. It
would be interesting to extend results in this paper to the general case.

\subsection{Coherence equations}
We provide an equivalent condition of the compatibility and coherency of unit
actions in terms of linear relations among the binary operations in an operad.
This criterion will be applied in Theorem~\mref{thm:crel} to classify
\Loday operads with coherent or compatible unit actions.
\begin{theorem}
Let $\calp=(\dfgen,\dfrel,\star)$ be an \Loday operad.
\begin{enumerate}
\item
A unit action $(\alpha,\beta)$ on $\calp$ is coherent
if and only if, for every
$$(\sum_{i}\dfoa_i\ot \dfob_i, \sum_j\dfoc_j\ot \dfod_j)\in \dfrel,$$
the following {\bf coherence equations} hold.
\begin{enumerate}
\item[\bf{(C1)}] $\sum_i \beta(\dfoa_i)\dfob_i = \sum_j \beta(\dfoc_j)\dfod_j,$
\item[\bf{(C2)}] $\sum_i \alpha(\dfoa_i)\dfob_i  =\sum_j \beta(\dfod_j) \dfoc_j, $
\item[\bf{(C3)}] $\sum_i \alpha(\dfob_i)\dfoa_i =\sum_j \alpha(\dfod_j) \dfoc_j,$
\item[\bf{(C4)}] $\sum_i \beta(\dfob_i)\dfoa_i=\sum_j \beta(\dfoc_j)\beta(\dfod_j) \star, $
\item[\bf{(C5)}] $\sum_i \alpha(\dfoa_i)\alpha(\dfob_i) \star =
    \sum_j \alpha(\dfoc_j) \dfod_j.$
\end{enumerate}
\item
A unit action $(\alpha,\beta)$ on $\calp$ is compatible
if and only if, for every
$$(\sum_{i}\dfoa_i\ot \dfob_i, \sum_{j}\dfoc_j\ot
\dfod_j)\in \dfrel,$$
equations (C1), (C2) and (C3) above hold.
\end{enumerate}
\mlabel{thm:coherent}
\end{theorem}
Before proving Theorem~\ref{thm:coherent}, we give examples to show how it
can be used to determine compatible and coherent unit actions.
\begin{exam} {\rm
Consider the dendriform dialgebra $\calp_D=(\dfgenb_D,\dfrelb_D)$ in Eq.(\ref{eq:dia}).
So $\dfgenb_D=\{\prec,\succ\}$ and
$$
\dfrelb_D=\{(\prec\otimes\prec,\prec\otimes(\prec+\succ)),
(\succ\otimes\prec,\succ\otimes\prec),
((\prec+\succ)\otimes\succ,\succ\otimes\succ)\}.
$$
Suppose $(\alpha,\beta)$ is a coherent unit action of $\calp_D$.
Then the three equations in $\dfrelb_D$ satisfy (C1) -- (C5).
Applying (C1) to the first relation in $\dfrelb_D$, we obtain $
\beta(\prec) \prec =\beta(\prec) (\prec +\succ)$. So $\beta(\prec)
\succ =0$. Therefore $\beta(\prec)=0$. Applying (C2) to the first
relation, we have $\alpha(\prec) \prec = \prec$ since
$\beta(\prec+\succ)=\beta(\star)=1$. Thus $\alpha(\prec)=1$.
Similarly, applying (C2) to the second equation, we have
$\alpha(\succ) \prec =\beta (\prec) \succ$. So $\alpha(\succ)=0$.
Applying (C1) to the third equation gives $\beta(\prec+\succ)
\succ =\succ = \beta(\succ) \succ$. Thus $\beta(\succ)=1$.
Therefore, the only coherent unit action on $\calp_D$ is the one
given in~\cite{Lo4}:
$$ \alpha(\prec)=\beta(\succ)=1,
\alpha(\succ)=\beta(\prec))=0.$$ Note that we have only used
(C1) -- (C3). So the above is also the only compatible unit action of
$\calp_D$. }
\end{exam}
We will comment on the trialgebra case in Corollary~\mref{co:three}.

\begin{exam} {\rm
We consider the 2-associative algebra in \cite{L-R3} and \cite{Pi}.
It is given by generators $\dfgenb=\{\ast, \cdot\}$ and relations
$$\dfrelb=\{ (\ast\ot \ast, \ast\ot \ast), (\cdot\ot \cdot, \cdot\ot \cdot)\}.$$
Consider the unit action $(\alpha,\beta)$ in \cite{Lo6} given by
$\alpha(\ast)=\alpha(\cdot)=\beta(\ast)=\beta(\cdot)=1$.
We show that the action is not coherent regardless of the choice of
the associative operation $\star$.
Suppose $(\alpha,\beta)$ is coherent. Then applying (C4) to
$(\ast \ot \ast,\ast \ot \ast)$, we have
$\beta(\ast) \ast = \beta(\ast)\beta(\ast) \star$. So $\ast =\star$.
Applying (C4) to $(\cdot\ot \cdot,\cdot,\ot \cdot)$, we have
$\beta(\cdot)\cdot=\beta(\cdot)\beta(\cdot) \star$. So $\cdot=\star.$
This is impossible. So
Theorem~\ref{thm:Loday} cannot be applied to give a Hopf algebra
structure on free 2-associative algebras.

However, by verifying (C1) - (C3), we see that the unit action is
compatible when $\star$ is taken to be $\ast$. Loday and
Ronco~\cite{L-R3,L-R4} have equipped a free 2-associative algebra
with a Hopf algebra structure with $\ast$ as the product.
This suggests a possible connection between
compatibility and Hopf algebras. } \mlabel{ex:twoass}
\end{exam}
\smallskip

Similar arguments show that the associative dialgebra
(see Example~\mref{ex:ad}) and the operads $\mathcal{X}_\pm$ in~\mcite{Lo7}
have no coherent unit actions.

\proofbegin
(1)\ Let $(\alpha,\beta)$ be a unit action on $\calp$. We say that
the unit action $(\alpha,\beta)$ is coherent with a relation
$(\sum_i \dfoa_i\otimes \dfob_i,\sum_j \dfoc_j\ot \dfod_j)$ in $\dfrel$
if, for any
$\calp$-algebras $A$ and $B$, the operations $\dfop\in \dfgen$,
when extended to $A\boxtimes B$ by Eq. (\ref{eq:tensorext}), still
satisfy the same relation. Then to prove the theorem, we only need
to prove that $(\alpha,\beta)$ is coherent with a given relation
$(\sum_i \dfoa_i\otimes \dfob_i,\sum_j \dfoc_j\ot \dfod_j)\in \dfrel$ if
and only if (C1) -- (C5) hold for this relation.

Further, by definition, $(\alpha,\beta)$ is coherent with
$(\sum_i \dfoa_i\otimes \dfob_i,\sum_j \dfoc_j\ot \dfod_j)$ means that,
for any $\calp$-algebra $A$ and $B$ and for
any $a,a',a''\in A_+$ and $b,b',b''\in B_+$ such that
at least one of $b,b',b''$ is not 1, we have the equation
\begin{equation}
\sum_i \big ((a\ot b) \dfoa_i (a'\ot b')\big) \dfob_i (a''\ot b'')
    = \sum_j (a\ot b) \dfoc_j \big ( (a'\ot b')\dfod_j (a''\ot b'')\big ).
\mlabel{eq:coh}
\end{equation}
Thus there are 7 mutually disjoint cases for the choice of such
$b,b',b''$: the case when none of $b,b',b''$ is 1, the three cases
when exactly one of $b,b',b''$ is 1, and the three cases when
exactly two of $b,b',b''$ are 1. Note that when none of $b,b',b''$
is 1, Eq. (\ref{eq:coh}) just means that $(\sum_i \dfoa_i\otimes
\dfob_i,\sum_j \dfoc_j\ot \dfod_j)$ is a relation for $\calp$, so is
automatic true. Thus to prove the theorem we only need to prove

\begin{enumerate}
\item[] Case 1.  Eq. (\ref{eq:coh}) holds for $b=1$, $b'\neq 1\neq b''$ if and only
if (C1) is true;
\item[] Case 2. Eq. (\ref{eq:coh}) holds for $b'=1$, $b\neq 1\neq b''$ if and only
if (C2) is true;
\item[] Case 3. Eq. (\ref{eq:coh}) holds for $b''=1$, $b\neq 1\neq b'$ if and only
if (C3) is true;
\item[] Case 4. Eq. (\ref{eq:coh}) holds for $b=b'=1$, $b''\neq 1$ if and only
if (C4) is true;
\item[] Case 5. Eq. (\ref{eq:coh}) holds for $b'=b''=1$, $b\neq 1$ if and only
if (C5) is true;
\item[] Case 6. Eq. (\ref{eq:coh}) holds for $b=b''=1$, $b'\neq 1$ if
(C1) is true.
\end{enumerate}

We first consider the three cases when exactly one of $b,b',b''$ is 1.
Then by the definition of the operation $\dfop$ on $A\boxtimes B$
in Eq. (\ref{eq:tensorext}), we can rewrite Eq. (\ref{eq:coh}) as
\begin{equation}
 \sum_i \Big (\big ( (a\tstar a' )\tstar a''\big ) \otimes \big ((b\tdfoa_i b') \tdfob_i
    b''\big)\Big) =
    \sum_j \Big (\big (a\tstar(a'\tstar a'')\big) \otimes \big (b\tdfoc_j (b'\tdfod_j
    b'' ) \big)\Big),
\mlabel{eq:toper}
\end{equation}
for all $1\ot 1\ot 1 \neq b\ot b'\ot b''\in (B_+)^{\ot 3}$.
Since $\tstar$ is associative, by the arbitrariness of $A$ and $a,a',a''\in A_+$
(say by taking $a=a'=a''=1$),
we see that Eq. (\ref{eq:toper}), and hence Eq. (\ref{eq:coh}), is equivalent to
\begin{equation} \sum_i \big ((b\tdfoa_i b') \tdfob_i
    b''\big) =
    \sum_j \big (b\tdfoc_j (b'\tdfod_j b'' ) \big).
\mlabel{eq:oneb}
\end{equation}

{\bf Case 1.} Assume $b=1$ and $b',b''\neq 1$.
Then Eq. (\ref{eq:oneb}) is
$$ \sum_i (1\tdfoa_i b') \tdfob_i
    b'' =
    \sum_j 1\tdfoc_j (b'\tdfod_j b'' ) $$
and, by Eq. (\ref{eq:ext}), this means
$$ \sum_i \beta(\dfoa_i) b'\dfob_i b'' =
    \sum_j \beta(\dfoc_j) (b'\dfod_j b'' ).$$
That is,
$$  b' \big (\sum_i \beta(\dfoa_i)\dfob_i\big) b'' =
     b' \big (\sum_j \beta(\dfoc_j)\dfod_j \big) b'',$$
for every $\calp$-algebra $B$ and $b',b''\in B$. By the following
Lemma~\mref{lem:elt}, this is true if and only if
$$ \sum_i \beta(\dfoa_i)\dfob_i=\sum_j \beta(\dfoc_j)\dfod_j,$$
giving (C1).

\begin{lemma}
For $\dfop_1,\dfop_2\in \dfgen$, we have $\dfop_1=\dfop_2$ if and only
if $a \dfop_1 a'=a \dfop_2 a'$ for all $\calp$-algebras $A$ and
$a,a'\in A$. \mlabel{lem:elt}
\end{lemma}
\proofbegin
The only if part is clear. Now suppose $a \dfop_1 a'=a \dfop_2 a'$
for all $\calp$-algebras $A$ and $a,a'\in A$. Let $A$ be the free
$\calp$-algebra on one generator $x$. Then
$$A=\bfk x \oplus \calp_2 \oplus \cdots .$$
Here $\calp_2=\oplus_{i} \bfk \dfgene_i$ in which $\{\dfgene_i\}$ is a basis
of $\dfgen$. Also, a binary operation $\dfop\in \dfgen$ acts on $A\ot A$
by $$ x \dfop x = \dfop\in \calp_2.$$
Thus we have
$$ \dfop_1 = x \dfop_1 x = x \dfop_2 x = \dfop_2.$$
This proves the if part.
\proofend

{\bf Case 2.} Assume $b'=1$ and $b,b''\neq 1$. As in Case 1, we have
\begin{eqnarray*}
{\rm Eq. (\ref{eq:oneb})} && \olrarrow{\ }
 \sum_i (b\tdfoa_i 1) \tdfob_i  b'' =
    \sum_j b\tdfoc_j (1\tdfod_j b'' ),\ \forall\ b,b''\in B,\\
&& \olrarrow{\rm (\ref{eq:ext})}
 \sum_i b\alpha(\dfoa_i) \dfob_i b'' =
    \sum_j b\beta(\dfod_j) \dfoc_j b'',\ \forall\ b,b''\in B,\\
&& \olrarrow{ } \sum_i \alpha(\dfoa_i) \dfob_i =
    \sum_j \beta(\dfod_j) \dfoc_j \\
&&  \olrarrow{ } {\rm (C2)}.
\end{eqnarray*}

{\bf Case 3.} Assume $b,b'\neq 1$ and $b''=1$. As in Case 1, we have
\begin{eqnarray*}
{\rm Eq. (\ref{eq:oneb})} && \olrarrow{\rm (\ref{eq:ext})}
 \sum_i (b\dfoa_i b') \alpha(\dfob) = \sum_j b\dfoc_j b' \alpha(\dfod_j),\
    \forall\ b,b'\in B, \\
&& \olrarrow{ } {\rm (C3)}.
\end{eqnarray*}

We next consider the three cases when exactly two of $b,b',b''$ are the identity.
\smallskip

{\bf Case 4.} Assume $b=b'=1$ and $b''\neq 1$.
Now Eq. (\ref{eq:toper}) does not apply. Directly from Eq.~(\ref{eq:tensorext}),
we see that Eq.~(\ref{eq:coh}) means
$$\sum_i \big ((a\tdfoa_i a') \tstar a'' \big )\otimes (1\tdfob_i b'')
=\sum_j \big ( a\tstar (a'\tstar a'')\big ) \otimes \big (1
\tdfoc_j (1\tdfod_j b'')\big ).$$ Then by Eq. (\ref{eq:ext}), we
equivalently have
$$\sum_i \big ((a\tdfoa_i a') \tstar a'' \big )\otimes (\beta(\dfob_i) b'')
=\sum_j \big ( a\tstar (a'\tstar a'')\big )
        \otimes (\beta(\dfoc_j)\beta(\dfod_j) b'').$$
This means, by moving the scalars $\beta(\dfop_i^{(j)})$
across the tensor product,
$$\sum_i \beta(\dfob_i)\big ((a\tdfoa_i a') \tstar a'' \big )\otimes b''
=\sum_j \beta(\dfoc_j)\beta(\dfod_j)\big ( a\tstar (a'\tstar a'')\big )
        \otimes  b''.$$
Since this is true for all $B$ and $b''\in B$, we equivalently have
\begin{equation}
\sum_i \beta(\dfob_i)\big ((a\tdfoa_i a') \tstar a'' \big )
=\sum_j \beta(\dfoc_j)\beta(\dfod_j)\big ( a\tstar (a'\tstar a'')\big ).
\mlabel{eq:c41}
\end{equation}
Taking $a''=1$, we have
\begin{equation}
\sum_i \beta(\dfob_i) (a\tdfoa_i a')
=\sum_j \beta(\dfoc_j)\beta(\dfod_j)  (a\tstar a').
\mlabel{eq:c42}
\end{equation}
Conversely, right multiplying $a''$ (by $\tstar$) to this equation and
using the associativity of $\tstar$, we obtain Eq.~(\ref{eq:c41}).
So Eq.~(\ref{eq:c41}) and Eq.~(\ref{eq:c42}) are equivalent.

When $a$ and $a'$ are in $A$, Eq. (11) becomes
$$
\sum_i \beta(\dfob_i) (a\dfoa_i a')
=\sum_j \beta(\dfoc_j)\beta(\dfod_j)  (a\star a'),
$$
or, equivalently,
$$
a \big (\sum_i \beta(\dfob_i) \dfoa_i \big)a'
=a\big ( \sum_j \beta(\dfoc_j)\beta(\dfod_j)\, \star \big) a'.
$$
By Lemma~\mref{lem:elt}, we have
\begin{equation}
\sum_i \beta(\dfob_i) \dfoa_i =\sum_j \beta(\dfoc_j)\beta(\dfod_j) \star.
\mlabel{eq:c43}
\end{equation}
This is (C4). To go backwards, assuming (C4), then by the linearity of the map
$\dfop \mapsto \tdfop$ in Eq.~(\ref{eq:ext}), we get Eq.~(\ref{eq:c42})
for all augmented $\calp$-algebras $A_+$.
So we are done with Case 4.
\smallskip

{\bf Case 5.} Assume $b\neq 1$ and $b'=b''=1$. This case is similar to Case 4.
By Eq. (\ref{eq:tensorext}) we have
$$\sum_i \big ( (a\tstar a')\tstar a''\big )
    \otimes \big((b \tdfoa_i 1)\tdfob_i 1\big )
=\sum_j \big (a\tstar (a' \tdfod_j a'') \big )\otimes (b \tdfoc_j 1).$$
Then by Eq (\ref{eq:ext}), we have
$$\sum_i \big ( (a\tstar a')\tstar a''\big )
    \otimes \big (b \alpha(\dfoa_i)\alpha(\dfob_i) \big )
=\sum_j \big (a\tstar (a' \tdfod_j a'') \big )\otimes (b
\alpha(\dfoc_j)).$$ Taking $a=1$, $a'\neq 1\neq a''$, and moving
the scalars $\alpha(\dfop_i^{(j)})$ across the tensor product, we
have
$$\sum_i \alpha(\dfoa_i) \alpha(\dfob_i)  (a'\star a'') \otimes b
= \sum_j \alpha(\dfoc_j) (a\dfod_j a') \otimes b. $$ Since this is
true for any $B$ and $b\in B$, we have
$$\sum_i \alpha(\dfoa_i)\alpha(\dfob_i)  (a'\star a'')
= \sum_j \alpha(\dfoc_j) (a\dfod_j a'), $$ and by the
arbitrariness of $A$ and $a',a''\in A$, we have
$$\sum_i \alpha(\dfoa_i)\alpha(\dfob_i) \star
= \sum_j \alpha(\dfoc_j) \dfod_j. $$
This is (C5). As in Case 4, all implications here can be reversed.

{\bf Case 6.} Assume $b=b''=1$ and $b'\neq 1$. Then Eq. (\ref{eq:toper}) still applies
and we get
$$ \sum_i (1\tdfoa_i b')\tdfob_i 1 =\sum_j 1\tdfoc_j (b' \tdfod_j 1)$$
and, by Eq. (\ref{eq:ext}), we get
$$ \sum_i \big (\beta (\dfoa_i) \alpha(\dfob_i)\big ) b'=
    \sum_j \big (\beta(\dfoc_j) \alpha(\dfod_j)\big ) b'.$$
Thus
$$ \sum_i \beta (\dfoa_i) \alpha(\dfob_i)=\sum_j \beta(\dfoc_j) \alpha(\dfod_j).$$
We note that this follows from applying $\alpha$ to (C1). So we do not get a new
relation.

\medskip

This completes the proof of (1) of Theorem~\ref{thm:coherent}.

\medskip

\noindent (2)\ We note that the precise meaning of compatibility
of the unit action $(\alpha,\beta)$ with a relation $
(\sum_i \dfoa_i\ot \dfob_i,\sum_j  \dfoc_j\ot \dfod_j)$ in the space of
relations $\dfrel$ of $\calp$ is the requirement that Eq.
(\ref{eq:oneb}) holds for any $\calp$-algebra $B$ in the above
proof. Thus the verification of the compatibility condition is equivalent to
the verification of (C1), (C2) and (C3). This proves (2) of
Theorem~\ref{thm:coherent}.
\proofend

\section{Operads with coherent unit actions}
\mlabel{sec:coh}
We now apply Theorem~\ref{thm:coherent} to classify \Loday operads with coherent
unit actions and compatible unit actions. We then discuss some special cases.

\subsection{The classifications}

We display the relations of \Loday operads $(\dfgen,\dfrel,\star)$
that admit a coherent or compatible unit action.

\begin{theorem} Let $\calp=(\dfgen,\dfrel,\star)$ be an \Loday operad of dimension $n$
(that is, $\dim\dfgen=n$).
\begin{enumerate}
\item There is a coherent unit action $(\alpha,\beta)$ on $\calp$ with
$\alpha\neq \beta$ if and only if there is
a basis $\{ \dfop_i\}$ of $\dfgen$ with $\star=\sum_i \dfop_i$ such that
$\dfrel$ is contained in the subspace $\dfrel_{n,\,\coh}'$ of
$\dfgen^{\ot 2} \oplus \dfgen^{\ot 2}$ with the basis
\begin{equation}
\dfrelb'_{n,\,\coh} := \left\{ \begin{array}{l}
(\star \ot \dfop_2, \dfop_2\ot \dfop_2), \\
 (\dfop_1 \ot \dfop_1, \dfop_1 \ot \star), \\
 (\dfop_i\ot \dfop_1,  \dfop_i\ot \dfop_1),\ 2\leq i\leq n, \\
 (\dfop_2 \ot \dfop_j, \dfop_2\ot \dfop_j),\ 3\leq j\leq n, \\
 (\dfop_1\ot \dfop_i, \dfop_i \ot \dfop_2),\ 3\leq i\leq n,\\
 (\dfop_i\ot \dfop_j, 0), \
 (0, \dfop_i \ot \dfop_j),\ 3\leq i,j\leq n.
\end{array} \right\}
\mlabel{eq:an=b}
\end{equation}
\mlabel{it:cohneq}
\item There is a coherent unit action $(\alpha,\beta)$ on $\calp$ with
$\alpha= \beta$ if and only if there is a basis $\{ \dfop_i\}$ of $\dfgen$
with $\star=\sum_i \dfop_i$ such that $\dfrel$
is contained in the subspace $\dfrel_{n,\,\coh}''$ of
$\dfgen^{\ot 2} \oplus \dfgen^{\ot 2}$ with the basis
\begin{equation}
\dfrelb_{n,\,\coh}'':= \left \{ \begin{array}{l}
(\dfop_1 \ot \star, \dfop_1\ot\star) +(\star\ot \dfop_1, \star \ot \dfop_1)
 - (\dfop_1 \ot \dfop_1, \dfop_1 \ot \dfop_1), \\
 (\dfop_i\ot \dfop_j, 0), \
 (0, \dfop_i \ot \dfop_j),\ 2\leq i,j\leq n.
\end{array} \right \}.
\mlabel{eq:a=b}
\end{equation}
\mlabel{it:coheq}
\item There is a compatible unit action $(\alpha,\beta)$ on $\calp$ with
$\alpha\neq \beta$ if and only if there is
a basis $\{ \dfop_i\}$ of $\dfgen$ with $\star=\sum_i \dfop_i$ such that $\dfrel$
is contained in the subspace $\dfrel'_{n,\,\comp}$ of
$\dfgen^{\ot 2} \oplus \dfgen^{\ot 2}$ with the basis
\begin{equation}
\dfrelb_{n,\,\comp}' := \left\{ \begin{array}{l}
( (\dfop_1+\dfop_2) \ot \dfop_2, \dfop_2\ot \dfop_2), \\
 (\dfop_1 \ot \dfop_1, \dfop_1 \ot (\dfop_1+\dfop_2)), \\
 (\dfop_i\ot \dfop_1,  \dfop_i\ot \dfop_1),\ 2\leq i\leq n, \\
 (\dfop_2 \ot \dfop_j, \dfop_2\ot \dfop_j),\ 3\leq j\leq n, \\
 (\dfop_1\ot \dfop_i, \dfop_i \ot \dfop_2),\ 3\leq i\leq n,\\
 (\dfop_i\ot \dfop_2, 0), (0, \dfop_1\ot \dfop_i),\ 3\leq i \leq n,\\
 (\dfop_i\ot \dfop_j, 0), \
 (0, \dfop_i \ot \dfop_j),\ 3\leq i,j\leq n.
\end{array} \right\}
\mlabel{eq:an=bm}
\end{equation}
\mlabel{it:compneq}
\item There is a compatible unit action $(\alpha,\beta)$ on $\calp$ with
$\alpha = \beta$ if and only if there is
a basis $\{ \dfop_i\}$ of $\dfgen$ with $\star=\sum_i \dfop_i$ such that $\dfrel$
is contained in the subspace
$\dfrel_{n,\,\comp}''$ of $\dfgen^{\ot 2} \oplus \dfgen^{\ot 2}$ with the basis
\begin{equation}
\dfrelb_{n,\,\comp}'':= \left \{ \begin{array}{l}
(\dfop_1\ot \dfop_1, \dfop_1\ot \dfop_1),\ \\
(\dfop_1 \ot \dfop_i, \dfop_1\ot \dfop_i)+(\dfop_i \ot
\dfop_1,\dfop_i\ot \dfop_1),\
2\leq i\leq n, \\
 (\dfop_i\ot \dfop_j, 0), \
 (0, \dfop_i \ot \dfop_j),\ 2\leq i,j\leq n.
\end{array} \right \}.
\mlabel{eq:a=bm}
\end{equation}
\mlabel{it:compeq}
\end{enumerate}
\mlabel{thm:crel}
\end{theorem}

\proofbegin
(\ref{it:cohneq})
Let $\calp=(\dfgen,\dfrel,\star)$ be an \Loday operad.
Suppose there is a basis
$\{ \dfop_i\}$ with $\star=\sum_i \dfop_i$ such that $\dfrel$ is contained
in the subspace of $\dfgen^{\ot 2} \oplus \dfgen^{\ot 2}$ generated by
the relations (\ref{eq:an=b}).
Define linear maps $\alpha, \beta: \dfgen \to \bfk$ by
$$ \alpha(\dfop_i)=\delta_{1,i} 1,\ \beta(\dfop_i)=\delta_{2,i}1,\ 1\leq i\leq n,$$
where $\delta_{i,j}$ is the Kronecker delta.
Then $\alpha(\star)=1 =\beta(\star)$ and $\alpha\neq \beta$.
It is straightforward to check that each element in (\ref{eq:an=b}) satisfies
the equations (C1)-(C5) in Theorem~\ref{thm:coherent}.
For example, we check the first element $(\star\ot \dfop_2,\dfop_2\ot \dfop_2)$
against (C1)-(C5) and see that
\begin{itemize}
\item
(C1) means $\beta(\star)\dfop_2=\beta(\dfop_2)\dfop_2$ which holds since
$\beta(\star)=1=\beta(\dfop_2)$;
\item
(C2) means $\alpha(\star)\dfop_2=\beta(\dfop_2)\dfop_2$ which holds since
$\alpha(\star)=1=\beta(\dfop_2)$;
\item
(C3) means $\alpha(\dfop_2)\star =\alpha(\dfop_2) \dfop_2$ which holds since
$\alpha(\dfop_2)=0$;
\item
(C4) means $\beta(\dfop_2)\star=\beta (\dfop_2)\beta(\dfop_2) \star$ which
holds since $\beta(\dfop_2)=1$;
\item
(C5) means $\alpha(\star)\alpha(\dfop_2)\alpha =\alpha(\dfop_2)\dfop_2$ which
holds since $\alpha(\dfop_2)=0$.
\end{itemize}
Therefore, each element in $\dfrel$  satisfies
the equations (C1)-(C5). Thus the unit action $(\alpha,\beta)$ is coherent with
$\calp$. This proves the ``if" part.

\medskip
To prove the ``only if" part, we assume that there is a unit action
$(\alpha,\beta)$ on $\calp$ that is coherent and $\alpha\neq \beta$.
In particular, $\alpha(\star)=\beta(\star)=1$.
Then there are direct sum decompositions
$$ \dfgen = \bfk\star \oplus \ker \alpha = \bfk \star \oplus \ker \beta.$$
This, together with $\alpha(\star)=\beta(\star)=1$, implies that
$\alpha = \beta$ if and only if $\ker \alpha =\ker \beta$.
So we have $\ker \alpha \neq \ker \beta$. In fact, $\ker\alpha \subsetneq \ker\beta$
and $\ker\beta \subsetneq \ker\alpha$ since $\dim\ker\alpha=\dim\ker\beta=n-1$.
So there are elements $\dfop_1\in \ker \beta$ such that $\alpha(\dfop_1)\neq 0$
and $\dfop_2\in \ker \alpha$ such that $ \beta(\dfop_2)\neq 0$.
By rescaling, we can assume that
$$ \alpha(\dfop_1)=1=\beta(\dfop_2),\ \alpha(\dfop_2)=0 =\beta(\dfop_1). $$

Note that $\ker(\alpha)$ and $\ker(\beta)$ are the solution spaces of
the linear equations
$$\dfop_1\alpha(\dfop_1) x_1+\cdots + \alpha(\dfop_n) x_n
(=\alpha(x_1 \dfop_1+\cdots+ x_n \dfop_n)) =0 $$
and
$$ \beta(\dfop_1)x_1 + \cdots + \beta(\dfop_n) x_n
(=\beta(x_1 \dfop_1+ \cdots x_n \dfop_n)) =0.$$
So $\ker \alpha \cap \ker \beta$ is the solution space of the linear system
$$
\left \{ \begin{array}{l}
\alpha(\dfop_1) x_1+\cdots + \alpha(\dfop_n) x_n =0,\\
\beta(\dfop_1)x_1 + \cdots + \beta(\dfop_n) x_n =0.
\end{array} \right .
$$
Since $\alpha\neq \beta$, we see that
 $\ker \alpha\cap \ker \beta$ has dimension $n-2$.
The intersection also contains
$\star-\dfop_1-\dfop_2$. So there is a basis $\dfop_3,\cdots, \dfop_n$ of
$\ker\alpha\cap\ker\beta$ such that
$$\star-\dfop_1-\dfop_2 = \dfop_3+\cdots + \dfop_n.$$
Therefore
$$ \star=\dfop_1+\cdots +\dfop_n.$$

Now any element $\dfrele\in \dfrel$ is of the form
\begin{equation}
\dfrele= (\sum_{i,j=1}^n a_{ij} \dfop_i \otimes \dfop_j,
\sum_{i,j=1}^n b_{ij}\dfop_i\otimes \dfop_j).
\mlabel{eq:lrel}
\end{equation}
Since the unit action $(\alpha,\beta)$ is coherent on $\calp$,
$\lambda$ satisfies each of the five coherence
equations in Theorem~\ref{thm:coherent}. For (C1), the equation is
$$ \sum_{ij} a_{ij} \beta(\dfop_i)  \dfop_j
        = \sum_{ij} b_{ij} \beta(\dfop_i)  \dfop_j.$$
By our choice of the basis $\{\dfop_i\}$, this means
$$ \sum_j a_{2,j} \dfop_j = \sum_j b_{2,j} \dfop_j.$$
Thus we have
\begin{equation}
 a_{2,j}=b_{2,j},\ 1\leq j\leq n.
\mlabel{eq:co1}
\end{equation}

Similarly, from (C3), we obtain
\begin{equation}
 a_{i,1} = b_{i,1}.
\mlabel{eq:co3}
\end{equation}

Applying (C2), we obtain
$$ \sum_{ij} a_{ij} \alpha(\dfop_i) \dfop_j = \sum_{i,j} b_{ij}\beta(\dfop_j) \dfop_i.$$
This gives
\begin{equation}
    a_{1,i} =b_{i,2}.
\mlabel{eq:co2}
\end{equation}

Applying (C4), we have
$$ \sum_{ij} a_{ij}\beta(\dfop_j)\dfop_i =
    \sum_{ij} b_{ij} \beta(\dfop_i)\beta(\dfop_j) \star$$
which means
$$ \sum_i a_{i,2} \dfop_i = b_{2,2} \star.$$
If $b_{2,2}=0$, then since $\{\dfop_i\}$ is a basis, we have
$ a_{i,2}=0,\ 1\leq i\leq n.$
If $b_{2,2}\neq 0$, then again since $\{\dfop_i\}$ is a basis and
$\star =\sum_i \dfop_i$ by construction, we must have $a_{i,2}=b_{2,2}$.
Thus we always have
\begin{equation}
 a_{i,2}=b_{2,2},\ 1\leq i\leq n.
\mlabel{eq:co4}
\end{equation}

As with (C4), applying (C5) gives
\begin{equation}
 a_{1,1} = b_{1,j},\ 1\leq j\leq n.
\mlabel{eq:co5}
\end{equation}

Some of the relations above are duplicated. For examples,
$a_{1,1}=b_{1,1}$ is both in Eq. (\ref{eq:co3}) and Eq. (\ref{eq:co5}).
To avoid this we list (\ref{eq:co4}) and (\ref{eq:co5}) first and list the rest
only when needed. Note also that the above relations only involve coefficients
with at least one of the subscripts in $\{1,2\}$. This means that there are
no restrictions among the relations
$$\{(\dfop_i \ot \dfop_j,0),\ (0, \dfop_i \ot \dfop_j) |\, 3\leq i,j\leq n\}.$$
Then we see that $\dfrele$ in Eq. (\mref{eq:lrel}) is of the following linear combination of
linearly independent elements in $\dfgen^{\ot 2}\oplus \dfgen^{\ot 2}$.
\allowdisplaybreaks{
\begin{eqnarray*}
\lambda&=&b_{2,2} \big ( \sum_{i=1}^n \dfop_i\ot \dfop_2,
\dfop_2\ot \dfop_2\big)
+ a_{1,1} \big ( \dfop_1 \ot \dfop_1,  \dfop_1 \ot \sum_{j=1}^n \dfop_j\big)\\
&& + \sum_{i=2}^n a_{i,1} \big( \dfop_i\ot \dfop_1, \dfop_i\ot \dfop_1\big)
    + \sum_{j=3}^n a_{2,j} \big ( \dfop_2 \ot \dfop_j, \dfop_2\ot \dfop_j\big)\\
&& + \sum_{i=3}^n a_{1,i} \big ( \dfop_1\ot \dfop_i,  \dfop_i \ot \dfop_2\big)\\
&& + \sum_{i,j=3}^n a_{i,j} \big (\dfop_i\ot \dfop_j, 0\big)
   + \sum_{i,j=3}^n b_{i,j} \big (0, \dfop_i \ot \dfop_j\big).
\end{eqnarray*}
}
Recall that $\star =\sum_i \dfop_i$. We see that $\dfrele$ is in
$\dfrel_{n,\,\coh}'$ defined by Eq. (\ref{eq:an=b}).
So $\dfrel\subseteq \dfrel_{n,\,\coh}'$.

\medskip

(\mref{it:coheq})
The proof is similar to the last part.
To prove the ``if" direction, suppose there is a basis
$\{ \dfop_i\}$ with $\star=\sum_i \dfop_i$ such that $\dfrel$ is contained
in the subspace of $\dfgen^{\ot 2} \oplus \dfgen^{\ot 2}$ generated by
the relations (\ref{eq:a=b}).
Define linear maps $\alpha= \beta: \dfgen \to \bfk$ by
$$ \alpha(\dfop_i)=\delta_{1,i}1,\ 1\leq i\leq n.$$
Then $\alpha(\star)=1 =\beta(\star)$.
It is straightforward to check that
elements in (\ref{eq:a=b}) satisfy the equations (C1)-(C5) in
Theorem~\ref{thm:coherent}. 
Therefore, each element in $\dfrel$  satisfies the equations (C1)-(C5).
Thus the unit action $(\alpha,\beta)$ is coherent.

Now we consider the ``only if" direction. Let $(\alpha,\beta)$ be
a coherent unit action on $\calp$ with $\alpha=\beta$.
Then we have
$\dfgen=\bfk\star \oplus \ker \alpha$. Let $\{\dfop_2,\cdots,\dfop_n\}$
be a basis of $\ker\alpha=\ker\beta$
and define $\dfop_1=\star - (\dfop_2 +\cdots + \dfop_n)$. We have
$\star = \sum_{i} \dfop_i$, $\alpha(\dfop_1)=1$ and
$\alpha(\dfop_i)=0, 2\leq i\leq n.$

Let
\begin{equation}
\dfrele= (\sum_{i,j=1}^n a_{ij} \dfop_i \otimes \dfop_j,
\sum_{i,j=1}^n b_{ij}\dfop_i\otimes \dfop_j)
\mlabel{eq:lrel2}
\end{equation}
be in $\dfrel$. Since $(\alpha,\beta)$ is coherent on $\calp$,
applying Theorem~\mref{thm:coherent}, we have
\begin{itemize}
\item
by (C1), $\sum_j a_{1,j} \dfop_j = \sum_j b_{1,j} \dfop_j$,
so $a_{1,j}=b_{1,j};$
\item
by (C2), $\sum_j a_{1,j} \dfop_j = \sum_i b_{i,1} \dfop_i$,
so $a_{1,j} =b_{j,1};$
\item
by (C3), $\sum_i a_{i,1} \dfop_i =\sum_i b_{i,1} \dfop_i$,
so $a_{i,1}=b_{i,1};$
\item
by (C4), $\sum_i a_{i,1} \dfop_i = b_{1,1} \star$,
so $a_{i,1}=b_{1,1};$
\item
by (C5), $ a_{1,1} \star = \sum_j b_{1,j} \dfop_j$,
so $a_{1,1}=b_{1,j}$.
\end{itemize}
Therefore,
$$ a_{i,1}=a_{1,j}=b_{i,1}=b_{1,j},\ 1\leq i,j\leq n.$$
Thus
\begin{eqnarray*}
&& (\sum_{i=1 {\rm\ or\ }j=1} a_{ij} \dfop_i \otimes \dfop_j,
\sum_{i=1 {\rm\ or\ } j=1} b_{ij}\dfop_i\otimes \dfop_j)\\
&=& a_{1,1} (\sum_{i=1 {\rm\ or\ }j=1} \dfop_i \otimes \dfop_j,
\sum_{i=1 {\rm\ or\ } j=1} \dfop_i\otimes \dfop_j)\\
&=&
 a_{1,1} \Big( (\dfop_1 \ot \star,  \dfop_1\ot\star)
+((\star-\dfop_1)\ot \dfop_1,  (\star-\dfop_1) \ot \dfop_1) \Big )\\
&=& a_{1,1} \Big( (\dfop_1 \ot \star,  \dfop_1\ot\star)
    +(\star\ot \dfop_1, \star \ot \dfop_1)
 -  (\dfop_1 \ot \dfop_1, \dfop_1 \ot \dfop_1) \Big).
\end{eqnarray*}
On the other hand, the coherence equations impose no restriction on other
elements in $\dfgen^{\ot 2}\oplus \dfgen^{\ot 2}$.
Therefore $\dfrelb$ in Eq. (\ref{eq:lrel2}) is a linear combination of the
elements in Eq. (\mref{eq:a=b}).
So $\dfrel$ is a subspace of the subspace $\dfrel_{n,\,\coh}''$.

(\mref{it:compneq}) The proofs of part (3) of
Theorem~\ref{thm:crel} follow from a similar analysis as for part
(1). But we only need to consider (C1)-(C3) which give relations
(\ref{eq:co1}), (\ref{eq:co3}) and (\ref{eq:co2}).

(\mref{it:compeq})
Likewise, for the proof of part (4), we only consider (C1)-(C3) in the proof
of part (2). Grouping the resulting relations, we obtain Eq. (\ref{eq:a=bm}).
\proofend

\subsection{Special cases}
We now consider the cases where $\dfgen$ is of dimension 2 and  3.
Suppose $\alpha\neq \beta$.
Then from Theorem~\ref{thm:crel} we easily check that
$$\dfrelb_{2,\,\coh}' = \dfrelb_{2,\,\comp}'=\{ (\star \ot \dfop_2, \dfop_2\ot \dfop_2),
 (\dfop_1 \ot \dfop_1, \dfop_1 \ot \star),
 (\dfop_2\ot \dfop_1, \dfop_2\ot \dfop_1) \}.
$$

Replacing $\dfop_1$ by $\prec$ and replacing $\dfop_2$ by $\succ$, we obtain
the following improvement of Proposition 1.2 in~\cite{Lo6}.
\begin{coro}
Let $\calp=(\dfgen,\dfrel,\star)$ be an \Loday operad with $\dim \dfgen =2$.
The following statements are equivalent.
\begin{enumerate}
\item There is a coherent unit action $(\alpha,\beta)$ on
$\calp$ with $\alpha\neq \beta$;
\item There is a compatible unit action $(\alpha,\beta)$ on
$\calp$ with $\alpha\neq \beta$;
\item There is a basis $(\prec,\succ)$ of
$\dfgen$ with $\star=\prec+\succ$ such that $\dfrel$
is contained in the subspace of $\dfgen^{\ot 2} \oplus \dfgen^{\ot 2}$
with the basis
$$ \{ (\star \ot \succ, \succ\ot \succ),
 (\prec \ot \prec, \prec \ot \star),
 (\succ\ot \prec, \succ\ot \prec) \}.$$
\end{enumerate}
\end{coro}

Next suppose $\alpha=\beta$.
\begin{coro}
Let $\calp=(\dfgen,\dfrel,\star)$ be an \Loday operad with $\dim \dfgen =2$.
\begin{enumerate}
\item There is a coherent unit action $(\alpha,\beta)$ on
$\calp$ with $\alpha= \beta$ if and only
there is a basis $(\dfop_1,\dfop_2)$
of $\dfgen$ with $\star=\dfop_1+\dfop_2$ such that $\dfrel$
is contained in the subspace of $\dfgen^{\ot 2} \oplus \dfgen^{\ot 2}$ with
the basis
$$ \{ (\dfop_1\ot \star,\dfop_1\ot \star)+(\dfop_2\ot\dfop_1, \dfop_2\ot\dfop_1),
 (\dfop_2 \ot \dfop_2,0), (0, \dfop_2 \ot \dfop_2)\}.
$$
\item
There is a compatible unit action $(\alpha,\beta)$ on
$\calp$ with $\alpha= \beta$ if and only if
there is a basis $(\dfop_1,\dfop_2)$
of $\dfgen$ with $\star=\dfop_1+\dfop_2$ such that $\dfrel$ is contained in
the subspace of $\dfgen^{\ot 2} \oplus \dfgen^{\ot 2}$ with the basis
$$ \{ (\dfop_1 \ot \dfop_1, \dfop_1\ot \dfop_1),
 (\dfop_1 \ot \dfop_2,\dfop_1\ot \dfop_2)+ (\dfop_2\ot \dfop_1, \dfop_2 \ot \dfop_1),
(\dfop_2 \ot \dfop_2, 0), (0, \dfop_2\ot \dfop_2) \}.
$$
\end{enumerate}
\end{coro}

\smallskip

Now we consider the dimension three case. From Theorem~\ref{thm:crel} we get
\begin{coro}
Let $\calp=(\dfgen,\dfrel,\star)$ be an \Loday operad with $\dim \dfgen =3$.
There is a coherent unit action $(\alpha,\beta)$ on
$\calp$ with $\alpha\neq \beta$ if and only if
there is a basis $(\prec,\succ,\circ)$ of $\dfgen$ with
$\star=\prec+\succ+\circ$ such that $\dfrel$ is contained
in the subspace with the basis
\begin{eqnarray*}
\dfrelb_{3,\,\coh}'=& \{(\prec\otimes\prec,\prec\otimes\star),
(\succ\otimes\prec,\succ\otimes\prec),
(\star\otimes \succ,\succ\otimes\succ),
(\succ\otimes\circ,\succ\otimes\circ), \\
& (\prec\otimes\circ,\circ\otimes\succ),(\circ\otimes\prec,\circ\otimes\prec),
(\circ\otimes\circ, 0), (0, \circ\otimes\circ)\}.
\end{eqnarray*}
\mlabel{co:three}
\end{coro}
\proofbegin
It follows from Theorem~\ref{thm:crel} by replacing $\{\dfop_1,\dfop_2,\dfop_3\}$
with $\{\prec, \succ,\circ\}$.
\proofend

Clearly the relations in Eq. (\ref{eq:tri}) and (\ref{eq:NS}) of
the dendriform trialgebra and NS-algebra, respectively, are
contained in $\dfrelb_{3,\,\coh}'$, so have coherent unit actions,
as were shown in \cite{Le2,Lo6}.
We also note that when $n\geq 3$, compatibility does
not imply coherency for unit actions.

\section{Coherent unit actions on products and duals of operads}
\mlabel{sec:prod}
We briefly recall the concept of the black square product of \Loday
operads~\cite{E-G1,Lo7}, and show that coherent and compatible unit
actions are preserved by the black square product. We then recall the concept of
the dual operad and show that, other than a trivial case, coherence and compatibility
are not preserved by taking the duals.

\subsection{Products of operads}
For \Loday operads $(\dfgen_1,\dfrel_1,\star_1)$ and $(\dfgen_2,\dfrel_2,\star_2)$,
and for $\dfop^{(i)}\in \dfgen_i,\ i=1,2$, we use a column vector
$\dfpair{\dfop_1}{\dfop_2}$ to denote the tensor product
$\dfop_1\otimes \dfop_2\in \dfgen_1\otimes \dfgen_2$.
For
$f_i=(\dfoa_i\otimes\dfob_i, \dfoc_j\otimes\dfod_j)
\in \dfgen_i^{\otimes 2} \oplus \dfgen_i^{\otimes 2},\ i=1,2,$
define
$$\dfpair{f_1}{f_2}=
 \left (
\dfpair{\dfoa_1}{\dfoa_2}
\otimes\dfpair{\dfob_1}{\dfob_2},
\dfpair{\dfoc_1}{\dfoc_2}
\otimes\dfpair{\dfod_1}{\dfod_2}\right ) \in
(\dfgen_1\otimes \dfgen_2)^{\otimes 2}\oplus (\dfgen_1\otimes \dfgen_2)^{\otimes 2}.
$$
This extends by bilinearity to all
$f_i\in \dfgen_i^{\otimes 2} \oplus \dfgen_i^{\otimes 2},\ i=1,2.$
We define a subspace of $(\dfgen_1\otimes \dfgen_2)^{\otimes
2}\oplus (\dfgen_1\otimes \dfgen_2)^{\otimes 2}$ by
$$\dfrel_1\dftimes \dfrel_2=
\left \{ \dfpair{f_1}{f_2} \Big | f_i \in \dfrel_i,\ i=1,2 \right \}.
$$
So $\dfrel_1\dftimes \dfrel_2$ is a space of relations for the operator space
$\dfgen_1\otimes \dfgen_2$.

It is easy to see~\cite{E-G1} that the operad
$\calp:= (\dfgen_1\otimes \dfgen_2,\dfrel_1\dftimes \dfrel_2)$
with the operation $\dfpair{\star_1}{\star_2}$ is also an \Loday operad
\cite{E-G1}, called the {\bf black square} product of $\calp_1$
and $\calp_2$ and denoted $\calp_1\dftimes \calp_2$.

We recall the following results from~\cite{E-G1} for later
reference.

\begin{prop}
\begin{enumerate}
\item The quadri-algebra of Aguiar and Loday~\cite{A-L},
defined by four binary operations $\nearrow\,,\nwarrow\,,\searrow\,,\swarrow\,$
and 9 relations, is isomorphic to the black square product
$\calp_D \dftimes \calp_D$ where $\calp_D=(\dfgenb_D,\dfrelb_D,\star_D)$
is the dendriform dialgebra in Eq. (\ref{eq:dia}).
\item
The ennea-algebra of Leroux~\cite{Le1}, defined by 9 binary operations
$$ \nwarrow\,, \uparrow\,,\nearrow\,, \prec\,,\circ\,,\succ\,,\swarrow\,,
\downarrow\,, \searrow$$ and 49 relations, is isomorphic to the
black square product $\calp_T\dftimes \calp_T$ where
$\calp_T=(\dfgenb_T,\dfrelb_T,\star_T)$ is the
dendriform trialgebra in Eq. (\ref{eq:tri}). \item The
dendriform-Nijenhuis algebra~\cite{Le2}, equipped with 9 binary
operations
$$
\nearrow\,,\ \searrow\,,\ \swarrow\,,\ \nwarrow\,,\ \uparrow\,,\ \downarrow\,,\
\tilde{\prec}\,, \ \tilde{\succ}\,,\ \tilde{\bullet}
$$
and 28 relations, is isomorphic to the black square product
$\calp_T \dftimes \calp_N$ where
$\calp_T=(\dfgenb_T,\dfrelb_T,\star_T)$ and
$\calp_N=(\dfgenb_N,\dfrelb_N,\star_N)$ are the dendriform
trialgebra and the NS-algebra in Eq.
(\mref{eq:NS}), respectively.
\item The
octo-algebra~\cite{Le3}, defined using 8 operations
$$\nearrow_i,\nwarrow_i,\swarrow_i,\searrow_{\,i}\,,\ i=1,2$$
and 27 relations, is isomorphic to the third power
$\calp_D \dftimes \calp_D\dftimes \calp_D$ of the dendriform dialgebra
$\calp_D=(\dfgenb_D,\dfrelb_D,\star_D)$.
\end{enumerate}
\mlabel{pp:prod}
\end{prop}

\subsection{Unit actions on products}
We now use Theorem~\ref{thm:coherent} to show that coherent unit actions on \Loday
operads are preserved by the black square product and thus give rise to Hopf algebra
structures on the free objects of the product operad.

For $i=1,2$, let $\calp_i:=(\dfgen_i,\dfrel_i,\star_i)$ be an \Loday operad
and let $(\alpha_i,\beta_i)$ be a unit action on $\calp_i$.
Let $\calp=\calp_1\dftimes \calp_2$ and define

$$ \alpha:=
\alpha_1\otimes\alpha_2\ ({\rm resp.\ } \beta:=
\beta_1\otimes \beta_2) : \dfgen_1\otimes \dfgen_2 \to \bfk$$
by
$$ \alpha\left (\dfpair{\dfop_1}{\dfop_2}\right)
    =\alpha_1(\dfop_1)\alpha_2(\dfop_2)\ ({\rm resp.\ }
\beta\left (\dfpair{\dfop_1}{\dfop_2}\right)
    = \beta_1(\dfop_1)\beta_2(\dfop_2)).
$$
Then $(\alpha, \beta)$ defines
a unit action on $\calp$.

\begin{theorem}
Let $\calp_i:=(\dfgen_i,\dfrel_i,\star_i),\ i=1,2,$ be
\Loday operads with coherent unit actions $(\alpha_i,\beta_i)$.
Then the unit action $(\alpha_1\otimes \alpha_2, \beta_1\ot
\beta_2)$ on the \Loday operad
$\calp:=\calp_1\dftimes \calp_2:=
(\dfgen_1\otimes \dfgen_2, \dfrel_1\dftimes \dfrel_2,
\dfpair{\star_1}{\star_2})$ is also coherent.
Therefore, The augmented free $\calp$-algebra $\calp(V)_+$ on a $\bfk$-vector
space $V$ is a connected Hopf algebra. \mlabel{thm:Hopf}
\end{theorem}
It will be clear from the proof that, if $\calp_1$ and $\calp_2$ have compatible
unit actions, then so does their product.

\proofbegin
By Theorem~\ref{thm:coherent}, we just need to verify that the relations in
$\dfrel_1\dftimes \dfrel_2$ satisfy (C1)--(C5) with the unit action $(\alpha,\beta)$
on $\calp$. We recall that a relation in $\dfrel_1\dftimes \dfrel_2$ is of the form
$$\dfpair{f_1}{f_2}:=
 \left (\sum_{i_1,i_2}\dfpair{\dfoa_{i_1}}{\dfoa_{i_2}}\otimes
\dfpair{\dfob_{i_1}}{\dfob_{i_2}},
\sum_{j_1,j_2}\dfpair{\dfoc_{j_1}}{\dfoc_{j_2}}\otimes
\dfpair{\dfod_{j_1}}{\dfod_{j_2}} \right )$$ for
$(\sum_{i_1} \dfoa_{i_1}\otimes
\dfob_{i_1}, \sum_{j_1} \dfoc_{j_1}\ot \dfod_{j_1})\in \dfrel_1$ and
 $(\sum_{i_2} \dfoa_{i_2}\otimes \dfob_{i_2},
    \sum_{j_2} \dfoc_{j_2}\ot \dfod_{j_2})\in \dfrel_2$.

Then by bilinearity of the product $\dfpair{\dfop_1}{\dfop_2}$ and
the equation (C1) for $\calp_1$ and $\calp_2$, we have
\allowdisplaybreaks{\begin{eqnarray*} \lefteqn{ \sum_{i_1,i_2}
\beta\left (\dfpair{\dfoa_{i_1}}{\dfoa_{i_2}}\right )
\dfpair{\dfob_{i_1}}{\dfob_{i_2}}
=\sum_{i_1,i_2} \beta_1(\dfoa_{i_1})\beta_2(\dfoa_{i_2})
\dfpair{\dfob_{i_1}}{\dfob_{i_2}}} \\
&=& \sum_{i_1,i_2}  \dfpair{\beta_1(\dfoa_{i_1})
    \dfob_{i_1}}{\beta_2(\dfoa_{i_2})\dfob_{i_2}} \\
&=& \dfpair{\sum_{i_1} \beta_1(\dfoa_{i_1})\dfob_{i_1}}
    {\sum_{i_2} \beta_2(\dfoa_{i_2})\dfob_{i_2}} \\
&=& \dfpair{\sum_{j_1} \beta_1(\dfoc_{j_1})\dfod_{j_1}}
    {\sum_{j_2} \beta_2(\dfoc_{j_2})\dfod_{j_2}} \\
&=& \sum_{j_1,j_2}  \dfpair{\beta_1(\dfoc_{j_1})\dfod_{j_1}}
    {\beta_2(\dfoc_{j_2})\dfod_{j_2}} \\
&=& \sum_{j_1,j_2} \beta_1(\dfoc_{j_1})\beta_2(\dfoc_{j_2})
    \dfpair{\dfod_{j_1}}{\dfod_{j_2}} \\
&=& \sum_{j_1,j_2}  \beta\left (\dfpair{\dfoc_{j_1}}{\dfoc_{j_2}}\right )
    \dfpair{\dfod_{j_1}}{\dfod_{j_2}} .
\end{eqnarray*}}
This gives (C1) for the product operad $\calp$.
Conditions (C2) and (C3) can be verified in the same way.

For (C4), we have
\allowdisplaybreaks{\begin{eqnarray*} \lefteqn{
\sum_{i_1,i_2} \beta\left (\dfpair{\dfob_{i_1}}{\dfob_{i_2}}\right )
\dfpair{\dfoa_{i_1}}{\dfoa_{i_2}}
=\sum_{i_2,i_2} \beta_1(\dfob_{i_1})\beta_2(\dfob_{i_2})
    \dfpair{\dfoa_{i_1}}{\dfoa_{i_2}}} \\
&=& \dfpair{\sum_{i_1} \beta_1(\dfob_{i_1})\dfoa_{i_1}}
    {\sum_{i_2} \beta_2(\dfob_{i_2})\dfoa_{i_2}} \\
&=& \dfpair{\sum_{j_1} \beta_1(\dfoc_{j_1})\beta_1(\dfod_{j_1})\star_1}
{\sum_{j_2} \beta_2(\dfoc_{j_2})\beta_2(\dfod_{j_2})\star_2}\\
& = & \sum_{j_1,j_2}
\dfpair{\beta_1(\dfoc_{j_1})\beta(\dfod_{j_1})\star_1}
{\beta_2(\dfoc_{j_2})\beta(\dfod_{j_2})\star_2} \\
&=& \sum_{j_1,j_2} \beta_1(\dfoc_{j_1})\beta_1(\dfod_{j_1})
    \beta_2(\dfoc_{j_2})\beta_2(\dfod_{j_2})
\dfpair{\star_1}{\star_2} \\
&=& \sum_{j_1,j_2} \beta\left (\dfpair{\dfoc_{j_1}}{\dfoc_{j_2}}\right )
\beta\left (\dfpair{\dfod_{j_1}}{\dfod_{j_2}}\right )\dfpair{\star_1}{\star_2}.
\end{eqnarray*}}
This gives (C4) for the product operad $\calp$. The same works for (C5).
\proofend

By Proposition~\ref{pp:prod} and Theorem~\ref{thm:Hopf} we immediately obtain
the following results on Hopf algebras.
\begin{coro}
\begin{enumerate}
\item {\bf (Loday\cite{Lo6})} Augmented free quadri-algebras are Hopf algebras;
\item {\bf (Leroux\cite{Le1})} Augmented free ennea-algebras are Hopf algebras;
\item {\bf (Leroux\cite{Le2})} Augmented free Nijenhuis-dendriform algebras are Hopf algebras;
\item {\bf (Leroux\cite{Le3})} Augmented free octo-algebras are Hopf algebras;
\end{enumerate}
\mlabel{co:Hopf}
\end{coro}


\subsection{Duality of \Loday operads}
\mlabel{sec:dual}
We now recall the definition of the dual~\cite{E-G1,Lo4,Lo7} of an \Loday operad before we study
the relation between coherent unit actions and taking duals.

For an \Loday operad $\calp=(\dfgen,\dfrel,\star)$,
let $\check{\dfgen}:=\Hom(\dfgen,\bfk)$ be the dual space of $\dfgen$, giving the
natural pairing
$$ \langle\ ,\  \rangle_\dfgen: \dfgen \times \check{\dfgen} \to \bfk.$$
Then $\check{\dfgen}^{\otimes 2}$ is identified with the dual space of $\dfgen^{\otimes 2}$,
giving the natural pairing
$$\langle\ ,\  \rangle^{\otimes 2}: \dfgen^{\otimes 2} \times \check{\dfgen}^{\otimes 2} \to \bfk,
\langle x\otimes y, a\otimes b\rangle^{\otimes 2}
=\langle x, a \rangle_\dfgen\ \langle y, b\rangle_\dfgen.$$
We then define a pairing
\begin{equation}
\langle\ ,\  \rangle: (\dfgen^{\otimes 2} \oplus \dfgen^{\otimes 2}) \times
    (\check{\dfgen}^{\otimes 2} \oplus \check{\dfgen}^{\otimes 2}) \to \bfk
\mlabel{eq:dual}
\end{equation}
by
$$\langle (\alpha,\beta), (\gamma,\delta)\rangle = \langle \alpha,\gamma\rangle^{\otimes 2}
-\langle \beta, \delta\rangle^{\otimes 2}, \alpha,\beta\in \dfgen^{\otimes 2},\
\gamma,\delta\in\check{\dfgen}^{\otimes 2}.$$

We now define $\dfrel^{\bot}$ to be the annihilator of
$\dfrel\subseteq \dfgen^{\otimes 2} \oplus \dfgen^{\otimes 2}$ in $\check{\dfgen}^{\otimes 2}\oplus
\check{\dfgen}^{\otimes 2}$ under the pairing $\langle\ , \rangle$.
We call $\calp^{\,!}:=(\check{\dfgen},\dfrel^{\bot})$ the {\bf dual operad} of
$\calp=(\dfgen,\dfrel)$ which is the Koszul dual in our special case.
It follows from the definition that $(\calp^{\,!})^!=\calp$.

\begin{exam}{\rm
(Associative dialgebra~\cite[Proposition 8.3]{Lo4})
Let $(\dfgen_D, \dfrel_D)$ be the operad for the dendriform dialgebra.
Let $\{\dashv,\vdash\}\in \check{\dfgen}_D$ be the dual basis of $\{\prec, \succ\}$
(in this order). Then $\dfrel_D^\bot$ is generated by
$$ \{ (\dashv\otimes \dashv, \dashv\otimes \dashv),
(\vdash\otimes \vdash, \vdash\otimes \vdash),
(\dashv \otimes\dashv, \dashv \otimes \vdash),
(\vdash \otimes \dashv, \vdash \otimes \dashv),
(\dashv \otimes \vdash, \vdash \otimes \vdash)\}.
$$
It is called the associative dialgebra and is denoted
$\calp_{AD}=(\dfgen_{AD},\dfrel_{AD})$. It has two associative
operations $\dashv$ and $\vdash$, giving two
\Loday operads $(\dfgen_{AD},\dfrel_{AD},\dashv)$ and
$(\dfgen_{AD},\dfrel_{AD},\vdash)$. They are isomorphic under the
correspondence $\dashv\, \leftrightarrow\, \vdash$. }
\mlabel{ex:ad}
\end{exam}
Other examples of duals of \Loday operads can be found
in~\cite{E-G1,L-R2}.

\subsection{Unit actions on duals}
A natural question to ask is whether the dual of an \Loday operad
$\calp=(\dfgen,\dfrel,\star)$ with a coherent unit action still has a coherent
unit action. When $\dfgen$ has one generator, the answer is
positive. This is because the generator must be the associative
operation $\star$. Then $(\dfgen,\dfrel,\star)$ is the operad for
associative algebras, with the coherent unit action given by
$\alpha (\star)=\beta(\star)=1$. The dual operad is isomorphic to
the operad itself, so  again admits a coherent unit action. We now show that
this is the only case that coherent unit action is preserved by taking the
dual.

\begin{theorem}
Suppose that an \Loday operad $\calp=(\dfgen,\dfrel,\star)$
with $n=\dim \dfgen\geq 2$
has a compatible unit action. Let $\calp^!=(\check{\dfgen},\dfrel^\bot)$ be
the dual operad.
\begin{enumerate}
\item $\calp^!$ has an associative operation, so is an \Loday
operad. In fact it has $n$ linearly independent
associative operations. \mlabel{it:dcoh} \item $\calp^!$ does not
have a compatible unit action for any choice of associative
operation. \mlabel{it:dcomp}
\end{enumerate}
The same is true when compatible is replaced by coherent.
\mlabel{thm:dual}
\end{theorem}
For example the associative dialgebra operad $\calp_{AD}$ in
Example~\mref{ex:ad} does not have a coherent unit action
even though it has two linearly independent associative binary operations
$\vdash$ and $\dashv$.

\proofbegin
(\mref{it:dcoh})\ First assume that the compatible unit
action on the given \Loday operad $(\dfgen,\dfrel)$ satisfies
$\alpha\neq \beta$. Let $\{\dfop_1,\cdots,\dfop_n\}$ be the bases
of $\dfgen$ given in the proof of Theorem~\ref{thm:crel}.(\mref{it:compneq}). So we
have
$$\sum_i \dfop_i =\star,\ \alpha(\dfop_i)=\delta_{1,i},\ \beta(\dfop_i)=\delta_{2,i}.$$
Then by Theorem~\ref{thm:crel}, $\dfrel$ is a subspace of
the subspace $\dfrel'_{n,\,\comp}$ of $\dfgen^{\ot 2} \oplus \dfgen^{\ot 2}$ with a basis
\begin{equation}
\dfrelb'_{n,\,\comp} := \left\{ \begin{array}{l}
( (\dfop_1+\dfop_2) \ot \dfop_2, \dfop_2\ot \dfop_2), \\
 (\dfop_1 \ot \dfop_1, \dfop_1 \ot (\dfop_1+\dfop_2)), \\
 (\dfop_i\ot \dfop_1,  \dfop_i\ot \dfop_1),\ 2\leq i\leq n, \\
 (\dfop_2 \ot \dfop_j, \dfop_2\ot \dfop_j),\ 3\leq j\leq n, \\
 (\dfop_1\ot \dfop_i, \dfop_i \ot \dfop_2),\ 3\leq i\leq n,\\
 ((\dfop_i\ot \dfop_2, 0), (0, \dfop_1\ot \dfop_i), 3\leq i \leq n,\\
 (\dfop_i\ot \dfop_j, 0), \
 (0, \dfop_i \ot \dfop_j),\ 3\leq i,j\leq n.
\end{array} \right\}
\mlabel{eq:an=bm2}
\end{equation}
Let $\{\check{\dfop}_i\}$ be the dual basis of $\{\dfop_i\}$. So
$\check{\dfop}_i (\dfop_j)=\delta_{i,j},\ 1\leq i,j\leq n.$

Then the pairing (\mref{eq:dual})
between $\dfgen^{\ot 2}\oplus \dfgen^{\ot 2}$ and
$\check{\dfgen}^{\ot 2} \oplus \check{\dfgen}^{\ot 2}$ is given by
$$\langle (\dfop_i\otimes \dfop_j, \dfop_k\otimes \dfop_\ell),
(\check{\dfop}_s\otimes \check{\dfop}_t, \check{\dfop}_u\otimes \check{\dfop}_v) \rangle
= \delta_{i,s}\delta_{j,t}-\delta_{k,u}\delta_{\ell,v}.$$
Consider $x=(\check{\dfop}_1\ot \check{\dfop}_1,\check{\dfop}_1\ot \check{\dfop}_1)$.
Then the pairing between $x$ and the second element in $\dfrele'_{n,\,\comp}$ is $1-1=0$.
The pairing between $x$ and every other element in $\dfrele'_{n,\,\comp}$ is also 0 since
$\check{\dfop}_1\ot \check{\dfop}_1$ does not occur in the element.
Thus $x=(\check{\dfop}_1\ot \check{\dfop}_1,\check{\dfop}_1\ot \check{\dfop}_1)$
is in $\dfrel^{\prime}_{n,\,\,\comp}{}^\bot$ and hence
in $\dfrel^\bot$, the relation space of $\calp^!$. This shows that $\calp^!$ has
an associative operation. The same argument works for
$(\check{\dfop}_i\ot \check{\dfop}_i,\check{\dfop}_i\ot \check{\dfop}_i)$,
$2\leq i\leq n$, giving $n$ linearly independent associative operations.

\smallskip

Next assume that the compatible unit action on the given \Loday
operad $(\dfgen,\dfrel)$ satisfies $\alpha= \beta$. Let
$\{\dfop_1,\cdots,\dfop_n\}$, $n=\dim \dfgen$, be the
basis of $\dfgen$ given in the proof of
Theorem~\ref{thm:crel}.(4). Then as in the last case, we check
that $(\check{\dfop}_i\ot \check{\dfop}_i,\check{\dfop}_i\ot
\check{\dfop}_i)$, $1\leq i\leq n$, are in
$\dfrel''_{n,\,\comp}{}^\bot\subseteq \dfrel^\bot$, again showing
that $\calp^!$ is \mbox{\Loday.}

\medskip

(2)
Suppose that there is a choice of associative operation $\star$ in $\check{\dfgen}$
and a unit action $(\alpha,\beta)$ that is compatible with
$(\check{\dfgen},\dfrel^\bot, \star)$.
First assume that $\alpha\neq \beta$. We easily verify that the three relations
\begin{equation}
(\check{\dfop}_2\ot \check{\dfop}_1,\check{\dfop}_2\ot \check{\dfop}_1),
(\check{\dfop}_1\ot\check{\dfop}_1,\check{\dfop}_1\ot\check{\dfop}_1),
(\check{\dfop}_2\ot\check{\dfop}_2,\check{\dfop}_2\ot\check{\dfop}_2)
\mlabel{eq:rel2}
\end{equation}
are in $\dfrel'_{n,\,\comp}{}^\bot$ and hence in $\dfrel^\bot$. So they should satisfy
the compatible equations (C1)-(C3) in Theorem~\ref{thm:crel}.  Applying (C2) to
the first relation in Eq. (\ref{eq:rel2}),  we obtain
$\alpha (\check{\dfop}_2) \check{\dfop}_1 = \beta(\check{\dfop}_1)\check{\dfop}_2$.
So
\begin{equation}
\alpha(\check{\dfop}_2)=\beta(\check{\dfop}_1)=0.
\mlabel{eq:AD}
\end{equation}
Applying (C2) to the second relation in Eq. (\ref{eq:rel2}), we
have $\alpha(\check{\dfop}_1)=\beta(\check{\dfop}_1)$, yielding
$\alpha(\check{\dfop}_1)=0$ by Eq. (\ref{eq:AD}). Applying (C2) to
the third relation in Eq. (\ref{eq:rel2}), we have
$\alpha(\check{\dfop}_2)=\beta(\check{\dfop}_2)$, giving
$\beta(\check{\dfop}_2)=0$ by Eq. (\ref{eq:AD}). For $i\geq 2$, we
check that $(\check{\dfop}_i\ot\check{\dfop}_2,0)$ and
$(0,\check{\dfop}_1\ot \check{\dfop}_i)$ are in $\dfrel'_{n,\,
\comp}{}^\bot$ and hence in $\dfrel^\bot$, so satisfy (C1)-(C2).
Applying (C1) to the first relation gives
$\beta(\check{\dfop}_i)\check{\dfop}_2=0$. So
$\beta(\check{\dfop}_i)=0$. Applying (C3) to the second equation
gives $\alpha(\check{\dfop}_i)\check{\dfop}_1=0$. So
$\alpha(\check{\dfop}_i)=0$. Therefore $\alpha$ and $\beta$ are
identically zero. But this is impossible, since $\alpha(\star)$
and $\beta(\star)$ should be 1.

\smallskip

Next assume that $\alpha=\beta$. For each $i\geq 2$ (there is such an $i$
since $\dim\dfgen \geq 2$), we check that
$(\check{\dfop}_1\ot\check{\dfop}_i, \check{\dfop}_1\ot \check{\dfop}_i)$
and
$(\check{\dfop}_i\ot\check{\dfop}_1, \check{\dfop}_i\ot \check{\dfop}_1)$
are in $\dfrel''_{n,\,\comp}{}^\bot\subseteq \dfrel^\bot$. Applying (C2) to them, we get
$$ \alpha(\check{\dfop}_1)\check{\dfop}_i=\beta(\check{\dfop}_i)\check{\dfop}_1,
\quad \alpha(\check{\dfop}_i)\check{\dfop}_1=\beta(\check{\dfop}_1)\check{\dfop}_i.$$
So $\alpha(\check{\dfop}_1)=\beta(\check{\dfop}_i)=0$ and
$\alpha(\check{\dfop}_i)=\beta(\check{\dfop}_1)=0, i\geq 2.$
Thus $\alpha$ and $\beta$ are identically zero, giving a contradiction.

\medskip
Finally if $\calp$ has a coherent unit action, then it
automatically has a compatible unit action. So by part (\mref{it:dcoh}),
$\calp^!$ is an \Loday operad, but does not have a compatible unit
action. It therefore does not have a coherent unit action.
\proofend

\medskip

{\em Acknowledgements}: The first author would like to thank the
Ev. Studienwerk Villigst and the theory department of the
Physikalisches Institut, at Bonn University for generous support.
The second author is supported in part by grants from NSF and
the Research Council of the Rutgers University.
We thank J.-L. Loday for suggestions and encouragement.

\addcontentsline{toc}{section}{\numberline {}References}

\end{document}